\documentclass[microtype]{gtpart}

\usepackage{pinlabel}  
\usepackage[utf8]{inputenc}
\usepackage{stmaryrd}
\usepackage{mathtools}
\usepackage{cancel} 
\usepackage{calc} 
\usepackage{tikz} 
\usepackage{xparse}
\usepackage{lmodern}

\title{Kontsevich's Swiss Cheese Conjecture}

\author{Justin D. Thomas}
\givenname{Justin D.}
\surname{Thomas}
\address{Department of Mathematics\\
University of Notre Dame\\
255 Hurley Hall\\Notre Dame\\IN\\46556-4618
}
\email{jthoma20@nd.edu}
\urladdr{http://www.nd.edu/~jthoma20/}


\arxivreference{math.AT/1011.1635}
\arxivpassword{wv4kc}

\volumenumber{}
\issuenumber{}
\publicationyear{}
\papernumber{}
\startpage{}
\endpage{}
\doi{}
\MR{}
\Zbl{}
\received{}
\revised{}
\accepted{}
\published{}
\publishedonline{}
\proposed{}
\seconded{}
\corresponding{}
\editor{}
\version{}

\theoremstyle{plain} 
\newtheorem*{lemma*}{Lemma}
\newtheorem*{theorem*}{Theorem}
\newtheorem*{remark*}{Remark}
\swapnumbers
\newtheorem{theorem}{Theorem}[subsection]
\newtheorem{lemma}[theorem]{Lemma}
\newtheorem{corollary}[theorem]{Corollary}
\newtheorem{proposition}[theorem]{Proposition}

\theoremstyle{definition} 
\newtheorem{example}[theorem]{Example}
\newtheorem{definition}[theorem]{Definition}    
\newtheorem{notation}[theorem]{Notation}
\newtheorem{relations}[theorem]{Relations}
\newtheorem{remark}[theorem]{Remark}

\DeclareMathOperator*{\colim}{colim}

\DeclareMathOperator*{\hocolim}{hocolim}
\DeclareMathOperator{\map}{map}

\tikzset{math nodes/.style={execute at begin node=$,execute at end node=$}}
\tikzset{display math nodes/.style={execute at begin node=$\displaystyle, 
execute at end node=$}}
\tikzset{map name/.style={font=\scriptsize}}
\tikzset{numbered picture/.style={baseline=(current bounding box.center)}} 
\usetikzlibrary{matrix,arrows,calc,decorations.pathmorphing,trees}

\usetikzlibrary{positioning}

\newcommand{\calC}{\mathcal{C}}
\newcommand{\calD}{\mathcal{D}}

\newcommand{\calO}{\mathcal{O}}
\newcommand{\calP}{\mathcal{P}}
\newcommand{\calQ}{\mathcal{Q}}

\newcommand{\aut}{\mathrm{aut}}
\newcommand{\dom}{\mathrm{dom}\,}

\newcommand{\col}{\mathrm{col}}

\newcommand{\rt}{\mathrm{rt}}

\NewDocumentCommand{\Op}{gg}
  {%
  \IfNoValueTF{#2}%
	{\IfNoValueTF{#1}%
	  {\mathrm{Op}}%
	  {\mathrm{Op}_{#1}}}%
	{\mathrm{Op}_{#1}^{#2}}%
  }%

\newcommand{\hoch}{\mathrm{Hoch}}
\newcommand{\hh}{\mathrm{HH}}
\newcommand{\ch}{\mathrm{CH}}

\newcommand{\Alg}{\mathrm{Alg}}

\newcommand{\fun}{\mathrm{Fun}}
\newcommand{\enD}{\mathrm{End}}

\newcommand{\ind}{\mathrm{Ind}}
\newcommand{\modules}{\mathrm{Mod}}

\newcommand{\SymSeq}{\mathsf{sSeq}}
\newcommand{\h}{\mathsf{h}}
\newcommand{\e}{\mathsf{f}}
\newcommand{\wedone}{\mathsf{E}_{\mathsf{d-1}}}
\newcommand{\wed}{\mathsf{E}_{\mathsf d}}

\newcommand{\id}{\mathsf{id}}

\newcommand{\topl}{\mathsf{Top}}

\newcommand{\reals}{\mathbb{R}}

\newcommand{\edone}{E_{d-1}}

\newcommand{\Operads}{\Op{}{}}

\newcommand{\ul}{\underline}
\newcommand{\subc}{_{\calC}}
\newcommand{\inv}{^{-1}}

\newcommand{\abs}[1]{\left| #1 \right|}
\newcommand{\inthom}{\underline{\hom}}
\newcommand{\modsim}{/\!\!\sim\,}

\newcommand{\smallbullet}{{\scriptstyle\bullet}}
\newcommand{\W}{W}

\newcommand{\swcheese}{{SC}_{\!d}}
\newcommand{\swcheesebullet}{{SC}_{\!d,\smallbullet}}
\newcommand{\wswcheese}{\mathsf{SC}_{\mathsf{d}}}

\newcommand{\wcheeseh}{\mathsf{SC}_{\mathsf{d}}^{\mathsf{h}}}

\newcommand{\schone}{\mathsf{SC}_{\mathsf d}^{\h 1}}
\newcommand{\sch}{{\wswcheese^{\h}}}
\newcommand{\schonebullet}{\mathsf{SC}_{\mathsf d,\smallbullet}^{\h 1}}
\newcommand{\schbullet}{{\mathsf{SC}_{\mathsf d, \smallbullet}^{\h}}}

\newcommand{\swisscheese}{\mathsf{SC}}
\newcommand{\schinf}{{\mathsf{SC}_{\mathsf{d}}^{\h\infty}}}

\newcommand{\rfor}{\mathsf{For}}
\newcommand{\for}{\mathsf{For}}
\newcommand{\trees}{\mathsf{Trees}}

\newcommand{\outverts}{\mathrm{out}}
\newcommand{\exEdges}{\textsl{Edge}}

\newcommand{\spaces}{\calC}
\newcommand{\Coll}[1]{#1^{\Sigma}}
\newcommand{\Collleqone}[1]{{#1^{\Sigma_{\leq 1}}}}
\newcommand{\E}{\mathsf{E}}
\newcommand{\inverts}{{\mathrm{in}}}
\mathchardef\dash="2D
\newcommand{\indec}{\mathrm{indec}}

\begin{document}

\begin{abstract}    
  We prove a conjecture of Kontsevich which states that if $A$ is an $\edone$ 
  algebra then the Hochschild cochain object of $A$ is the universal $E_d$ 
  algebra acting on $A$.  The notion of an $E_d$ algebra acting on an $E_{d-1}$ 
  algebra was defined by Kontsevich using the swiss cheese operad of Voronov.  
  The degree 0 and 1 pieces of the swiss cheese operad can be used to build a 
  cofibrant model for $A$ as an $E_{d-1}-A$ module.  The theorem amounts to the 
  fact that the swiss cheese operad is generated up to homotopy by its degree 0 
  and 1 pieces.
\end{abstract}

\begin{asciiabstract}
  We prove a conjecture of Kontsevich which states that if A is an E(d-1) 
  algebra then the Hochschild cochain object of A is the universal E(d) algebra 
  acting on A.  The notion of an E(d) algebra acting on an E(d-1) algebra was 
  defined by Kontsevich using the swiss cheese operad of Voronov.  The degree 0 
  and 1 pieces of the swiss cheese operad can be used to build a cofibrant 
  model for $A$ as an $E_{d-1}-A$ module.  The theorem amounts to the fact that 
  the swiss cheese operad is generated up to homotopy by its degree 0 and 1 
  pieces.
\end{asciiabstract}

\maketitle

\section{Introduction}
\label{section:introduction}
In \cite{gerstenhaber-63} Gerstenhaber showed that the Hochschild cohomology 
$\hh^*(A)$ of an associative algebra $A$ is a graded Lie algebra and a graded 
commutative algebra, and the two structures are compatible.  Any graded vector 
space with this algebraic structure is now called a Gerstenhaber algebra.
In \cite{f-cohen-76} Cohen showed that the homology of the little disks operad, 
$H_{*}(E_2)$, is the Gerstenhaber operad.  Sinha also has shown this in 
\cite{sinha-06}.
Deligne later asked if the action of $H_{*}(E_2)$ on $\hh^*(A)$ descends from a 
natural action at the level of chains.  In other words, is there a natural 
algebra structure on $\ch^*(A)$ of $\mathsf{Chains}(E_2)$ which recovers the 
structure discovered by Gerstenhaber after passing to (co)homology?

Already, this question is evidently in the realm of homotopy theory.  So let us 
replace the associative algebras by $E_1$ algebras. This makes it clear that 
the question is fundamentally one about the relationship between the operads 
$E_1$ and $E_2$.  Indeed, we can generalize to consider the relationship 
between $E_d$ and $E_{d-1}$ algebras.
For any $\edone$ algebra in a sufficiently rich homotopical category $\spaces$ 
we can make sense of its Hochschild cochains as an object of $\spaces$. The 
Hochschild cochain object of $A$ is denoted $\hoch(A)$ and is an object of the 
same category to which $A$ belongs.  This terminology is based on the
case where $\spaces$ is the category of differential graded vector spaces and 
$A$ is an associative algebra.  In that case, $\hoch(A)$ is the usual 
Hochschild cochain complex of $A$.

The original Deligne conjecture where $A$ is an $E_1$ algebra in the category 
of chain complexes has been solved by Tamarkin \cite{tamarkin-operadic-proof}, 
Kontsevich-Soibleman \cite{kontsevich-soibelman-deformations}, 
Voronov \cite{voronov-homotopy-gerstenhaber}, 
McClure-Smith \cite{mcclure-smith}, 
Berger-Fresse \cite{bf-04-deligne}, 
and
Kaufmann-Schwell \cite{kaufmann-schwell-ainf-deligne}
.
A cyclic version is also due to Kaufmann \cite{kaufmann-cyclic}.
Vallette \cite{vallette-manin-products} generalized the theorem to include 
certain other Koszul operads.
The generalized version where $A$ is an $E_d$ algebra in a general category 
like $\spaces$ has been proven by Hu-Kriz-Voronov \cite{hu-kriz-voronov-2006} 
and in the $\infty$-operad setting by Lurie \cite{Lur2011}.  We show here that 
$\hoch(A)$ is not just an $E_d$ algebra, but comes equipped with a universal 
property.  It is the universal $E_d$ algebra acting on the $\edone$ algebra 
$A$.  This universal property is shown in the case $d = 2$ in a paper of 
Dolgushev-Tamarkin-Tsygan \cite{dolgushev-tamarkin-tsygan-2011}.

The notion of an $E_d$ algebra acting on an $\edone$ algebra was introduced by 
Kontsevich in \cite{kontsevich-operads-motives}.  This notion uses the swiss 
cheese operad $\swcheese$ of Voronov \cite{voronov-swiss-cheese}.  This is a 
two-colored operad which interpolates between $E_d$ and $\edone$.
A swiss cheese algebra is a pair $(B,A)$ where $B$ is an $E_d$ algebra, $A$ is 
an $\edone$ algebra, and there is some extra structure compatible with these as 
seen in definition~\ref{definition:ed-algebra-action-on-edone-algebra}.  We 
refer to this extra structure as an \emph{action} of $B$ on $A$.

The case $d=1$ is enlightening.  For simplicity, let us work in the category of 
vector spaces.  A (non-unital) $E_0$ algebra $A$ in vector spaces is just a 
vector space with no extra data.  The Hochschild object in this case is 
$\hom(A,A)$.  If $B$ is an associative algebra, it is in particular an $E_1$ 
algebra. An $SC_1$ structure on the pair $(B,A)$ then amounts to the choice of 
a $B$-module structure on $A$.

In this case, the swiss cheese conjecture merely states that $\hom(A,A)$ is an 
associative algebra, and giving  an $SC_1$ structure on $(B,A)$ is equivalent 
to giving a map of associative algebras $B \to \hoch(A) = \hom(A,A)$.
We prove the analog of this when $B$ is an $E_d$ algebra and $A$ is an $\edone$ 
algebra, $d \geq 1$.  

Joseph Hirsh brought to the attention of the author the following helpful 
characterization of the results of this paper. Given a bifibrant $\edone$ 
algebra $A$ the functor from $E_d$ algebras to spaces,
\[
	B \mapsto \{ \mbox{the space of swiss cheese actions of } B \mbox{ on } A\}
\]
is represented by the Hochschild cochain object of $A$.

\subsection{Acknowledgements} The author would like to thank Kevin Costello, 
Paul Goerss, John Francis, Mike Hopkins, Jacob Lurie, Vasily Dolgushev, Bill 
Dwyer, and Stephan Stolz for helpful conversations.  In addition, thanks go to 
Haynes Miller and  Ezra Getzler for many helpful remarks and to an anonymous 
referee for a large
number of useful comments which have greatly improved the presentation of 
results in this article.

\subsection{Outline of the paper}
\label{subsection:outline-of-the-paper}
In section~\ref{section:the-swiss-cheese-operad} we define the $E_d$ and 
$\swcheese$ operads, and give both an imprecise and precise statement of the 
theorem we will prove.  We also outline the idea of the proof.  In 
section~\ref{section:hochschild-cohomology-from-swiss-cheese} we define the 
Hochschild cochain object for $\edone$ algebras and show that we can use the 
swiss cheese operad to construct a model for the Hochschild cochain object.  We 
use this model in section~\ref{section:universal-cheese} to prove a ``universal 
cheese'' theorem which applies to an arbitrary operad acting on Hochschild 
cochains.  In section~\ref{section:right-Ed-action} we show that $E_d$ does 
indeed act, up to homotopy, on Hochschild cochains, and the ``universal 
cheese'' theorem specializes to the main theorem of this paper.  Finally, in 
section~\ref{section:the equivalence schinf to sch} we prove the main theorem 
which allows the homotopy $E_d$ action: the swiss cheese operad is freely 
generated up to homotopy by its degree 0 and 1 pieces.

\section{The swiss cheese operad}
\label{section:the-swiss-cheese-operad}
We will define $K$-colored operads in general, and the swiss cheese operad in 
particular.  We also describe algebras over the swiss cheese operad and state 
the main theorem of this paper.

\subsection{The colored operad swiss cheese}
	Fix a set $K$. A $K$-colored set is a pair $( I, \mathrm{col})$ where $I$ 
	is a set and $\mathrm{col}: I \rightarrow K $ is a map of sets, called the 
	\emph{coloring}.  We will often denote such a colored set simply by $I$, 
	leaving the coloring implicit.  Let $\aut ( I )$ be the group of bijections 
	on the set $I$ which preserve its coloring.  Since we have left the 
	coloring implicit, we use the notation $I^\#$ to refer to the underlying 
	uncolored set.

Let $( \mathcal S, \otimes)$ be a symmetric monoidal category.  We can speak of 
categories enriched over $\mathcal S$.  In particular, suppose $\mathcal O$ is 
a category enriched over $\mathcal S$ and suppose the objects of $\mathcal O$ 
are finite $K$-colored sets.  We let $\mathcal O(I;J)$ denote the $\mathcal S$ 
object of morphisms in $\mathcal O$ from $I$ to $J$.  If we further suppose 
that disjoint union of finite $K$-colored sets extends to an $\mathcal 
S$-enriched symmetric monoidal structure on $\mathcal O$, then each $\mathcal 
O(I;J)$ is a right $\aut(I)$ and left $\aut(J)$ module in a natural way.  In 
addition, if $\mathcal S$ contains finite coproducts, the symmetric monoidal 
structure on morphisms in $\mathcal O$ is specified by equivariant maps
\begin{equation} \label{equation:symmetric monoidal structure on enriched 
	morphisms}
	\coprod_{f\colon I^\# \to J^\#} \bigotimes_{j\in J} \mathcal O (f\inv\{j\};
	\{j\}) \to \mathcal O (I;J).
\end{equation}
The following definition follows Boardman-Vogt \cite{boardman-vogt-73}.  
\begin{definition}
\label{definition:colored-operad}
Let $( \mathcal S, \otimes)$ be a symmetric monoidal category containing all 
finite coproducts where $\otimes$ distributes over finite coproducts.  Let $K$ 
be a set.  The data of a \emph{$K$-colored operad $\mathcal O$ in the symmetric 
monoidal category $\mathcal S$} is a symmetric monoidal category, denoted 
$\mathcal O$, which is enriched over $\mathcal S$ and whose objects are 
$K$-colored finite sets, $I \rightarrow K$.  This data must satisfy the 
following conditions.  First, on objects, the symmetric monoidal structure of 
$\mathcal O$ is the disjoint union of sets over $K$.  Second, the map in 
\eqref{equation:symmetric monoidal structure on enriched morphisms} must be an 
isomorphism for every $I$ and $J$.
\end{definition}

\begin{remark}
We will use the unmodified noun \emph{operad} to mean $K$-colored operad when 
the coloring set $K$ is clear from context.  The reader should note that this 
differs from an equally plausible convention where \emph{operad} is always used 
to denote $\{*\}$-colored operads.
\end{remark}

\begin{notation}
\label{notation:underline-notation}
Let $\underline n$ to denote the finite set $\{1, \ldots, n \}$.
Typically $K$ will be $K=\{\e,\h\}$, where $\e$ stands for \emph{full disk} and 
$\h$ stands for \emph{half disk}.  In this case we use $(n,m)$
to denote the $K$-colored set which is the disjoint union of
  \begin{equation*}
\underline n \rightarrow \{\e\} \text{ and } \underline m \rightarrow\{\h\}.
\end{equation*}
If $K \simeq \{\ast\}$, then a $K$-colored operad will simply be called
a 1-colored operad.   Any 1-colored operad $\mathcal E$ gives for each $m, m' 
\in \mathbb Z_{\geq 0}$, objects $\mathcal E(\underline m'; \underline m) \in 
\mathcal S$.  We denote $\mathcal E(\ul m';\ul m)$ simply by $\mathcal E(m';m)$ 
and $\mathcal E( m; 1)$ simply by $\mathcal E(m)$.

Any $\{\e,\h\}$-colored operad $\mathcal O$ gives for each $n, m, n', m' \in 
\mathbb Z_{\geq 0}$, objects
\[
	\calO(n',m';n,m) \coloneqq \mathcal O((n',m'); (n, m)) \in \mathcal S.
\]
We denote $\mathcal O( n,m;0,1 )$ by $\mathcal O^\h ( n, m )$ and we denote 
$\mathcal O ( n,m;1,0 )$ by $\mathcal O^\e ( n, m )$.
\end{notation}


\begin{definition}
  Let $(\topl, \times)$ denote the symmetric monoidal category of compactly 
  generated Hausdorff topological spaces with the cartesian product.
\end{definition}

\begin{example}
\label{example:little-disks-operad}
Let $K$ be the one-point set $\{\e\}$ and let
 $(\mathcal S, \otimes) =(\topl, \times)$.   The operad $E_d$ is an 
 $\{\e\}$-colored operad in the category $\topl$.  Let $\bar D^d$ be the closed 
 unit disk inside $\mathbb R^d$. Call a map $f: \bar D^d \rightarrow \bar D^d$ 
 a \emph{little full-disk} (or \emph{little $d$ disk} or simply \emph{little 
 disk}) if $f$ is of the form $f ( x ) = r x + c$ for some $0 <r \leq 1$ and 
 $c\in \mathbb R^d$. Given a finite set $I$, the underlying set of $E_d ( I; \e 
 )$ is the set of embeddings
\[
f \colon \coprod_{i \in I} \bar D^d \to \bar D^d,
\]
where each restriction $f_i \colon \bar D^d \to \bar D^d$ is a little 
full-disk. Using notation~\ref{notation:underline-notation}, any isomorphism $I 
\to \underline n$ induces an isomorphism $E_d(I; \e) \simeq E_d(n)$, and the 
latter can naturally be considered as a subset of $\mathbb R^{n + d n}$.  This 
gives each $E_d(I; \e)$ a topology.  The operadic structure is given by 
composing little $d$ disks as maps $\bar D^d \rightarrow\bar D^d$.  The 
identity of $E_d$ is the little full-disk $\id \colon \bar D^d \to \bar D^d$.  
This is the unital version of $E_d$, so $E_d(0) = *$ and $E_d(1)$ consists of 
more than just the identity.
\end{example}

\begin{example}
Let $K = \{\e ,\h\}$. The $K$-colored operad $\swcheese$ is called the 
($d$-dimensional) swiss cheese operad and is the principal subject of this 
paper.  Like example~\ref{example:little-disks-operad} it is an operad in 
$(\topl,\times)$.

Let $\{\e\}$ and $\{\h\}$ denote the evident singleton $K$-colored sets.  By 
definition~\ref{definition:colored-operad} and 
formula~\eqref{equation:symmetric monoidal structure on enriched morphisms}, we 
only need to define the spaces $\swcheese(I;\e)$ and $\swcheese(I;\h)$ for 
every $K$-colored set $(I, \mathrm{col} \colon I \to K)$.  First, we define the 
``full-disk output'' part of $\swcheese$,
\[
 \swcheese ( I ; \e ) = \begin{cases}
						E_d(I; \e) & \mathrm{col}^{-1}(\e) = I \\
			\emptyset & \mbox{else}.
					   \end{cases}
\]
To define the ``half-disk output'' part of $\swcheese$, that is $\swcheese(I; 
\h)$, we first need the notion of \emph{little half-disks}.  Let $\bar D^d_+$ 
be the closed $d$-dimensional half-disk,
\[
\bar D^d_+ = \{ (x_1, \ldots, x_d) \in \mathbb R^d \mid \abs{x} \leq 1 \text{ 
and } x_d \geq 0 \}.
\]
A \emph{little half-disk} is defined to be a map $f \colon \bar D^d_+ \to \bar 
D^d_+$ of the form $f(x) = rx + c$ for some $0 < r \leq 1$ and $c \in 
\reals^{d-1} \times \{0\}$.
As a set, we define $\swcheese(I; \h)$ to consist of embeddings
\[
f \colon \coprod_{i \in I} \bar D_i^d \to \bar D^d_+
\]
where each restriction $f_i \colon \bar D_i^d \to \bar D^d_+$ is a little full 
disk (example \ref{example:little-disks-operad}) if $\mathrm{col}(i) = \e$ or a 
little half-disk if $\mathrm{col}(i) = \h$.  It is clear that if 
$\abs{\mathrm{col}^{-1}(\e)} = n$ and $\abs{\mathrm{col}^{-1}(\h)} = m$, then 
$\swcheese ( I; \h )$ can be naturally embedded inside $\mathbb R^{N}$ where $N 
= (d+1)n + dm $.  We give $\swcheese(I;\h)$ the subspace topology inherited 
from such an embedding.

Following notation~\ref{notation:underline-notation},
a point in $\swcheese^{\h}(n,m)$ is given by $n$ labeled full-disks and $m$ 
labeled half-disks in the unit half-disk where none of the disks intersect and 
the half-disks all lie on the bottom.  We allow the degenerate configuration 
when $(n,m) = (0,1)$ which is the unit half-disk contained in itself.  Note 
that we have $\swcheese^{\h}(0,0) = \ast$ and $\swcheese^{\h}(1,0)$ contains 
more than one point.  Thus we are using the unital swiss cheese operad.  This 
differs from Kontsevich  \cite{kontsevich-operads-motives} and Voronov 
\cite{voronov-swiss-cheese}.

Composition in $\swcheese$ is given by substituting full-disks and half-disks 
into each other 
.  More precisely, we have maps
\begin{equation}
 E_d(n) \times E_d(k_1) \times \cdots \times E_d(k_n) \to E_d(k_1 + \cdots + 
 k_n)
\end{equation}
and
\begin{multline}
\label{equation:swiss-cheese-full-composition}
 \swcheese^{\h}(n,m) \times
 E_d(k_1) \times \cdots \times E_d(k_n)  \times  \swcheese^{\h}(k_{n+1}, 
 \ell_1) \times \cdots \times \swcheese^{\h}(k_{n+m}, \ell_m)  \\
 \longrightarrow\swcheese^{\h}(k_1 + \cdots + k_{n+m}, \ell_1 + \cdots + 
 \ell_{m}).
\end{multline}


Notice that we can identify $\swcheese^{\h}(0,m)$ with $\edone(m)$ so that the 
restriction of $\swcheese$ to the spaces $\swcheese^{\h}(0,\smallbullet)$ is 
isomorphic to the operad $\edone$.  We say that $\edone$ is the $\h$ color of 
$\swcheese$ and $E_d$ is the $\e$ color of $\swcheese$.  We think of 
$\swcheese$ as interpolating between $E_d$ and $\edone$.
\end{example}

\begin{definition}
\label{definition:algebra-over-operad}
Suppose $\mathcal O$ is a $K$-colored operad in $\mathcal S$ and $\calC$ 
is a symmetric monoidal category enriched over $\mathcal S$.
An \emph{algebra} over $\mathcal O$ in the category $\calC$ is a strong 
symmetric monoidal functor $\mathcal O \rightarrow \calC$.  A morphism of 
$\mathcal O$ algebras is a monoidal natural transformation.  The category of 
$\mathcal O$ algebras in $\calC$ will be denoted $\Alg_{\mathcal 
O}(\calC)$.
\end{definition}

\begin{example}
\label{example:swiss-cheese-algebras}
If $\calC$ is a symmetric monoidal category enriched over $\topl$, we can 
consider  algebras over $\swcheese$ in $\calC$.   Such an algebra gives 
the data of a pair $(B,A)$ of objects in $\calC$ together with maps of 
topological spaces
\[
 E_d(n) \to \map\subc( B^{\otimes n}, B)
\]
and
\[
 \swcheese^{\h}(n,m) \to \map\subc( B^{\otimes n} \otimes A^{\otimes m}, A),
\]
where $\map\subc(C,C')$ is the topological space of maps between two objects 
$C, C'$ in $\calC$.
These data must satisfy conditions guaranteeing they assemble into a strong 
symmetric monoidal functor $\swcheese \to \calC$.

The object $B$ corresponds to the object $\{\e\}$ of $\swcheese$ and the object 
$A$ corresponds to $\{\h\}$.  Together these form a $K$-colored operad 
$\enD(B,A)$ in $\topl$ where, using  
notation~\ref{notation:underline-notation},
\begin{align}
\label{equation:endomorphism-operad-h}
 \enD^{\e}(B,A)(n,m) &= \map(B^{\otimes n} \otimes A^{\otimes m}, B) \\
 \label{equation:endomorphism-operad-e}
 \enD^{\h}(B,A)(n,m) &= \map(B^{\otimes n}\otimes A^{\otimes m}, A).
\end{align}
The collection of $\swcheese$ algebra structures on a fixed pair $(B,A)$ is the 
collection of strong symmetric monoidal functors
\[
 \swcheese \to \enD(B,A),
\]
which restrict to the identity on the set of objects.
Simply put, a swiss cheese algebra $(B,A)$ is an $E_d$ algebra $B$, an $\edone$ 
algebra $A$, and some chosen mixing of these structures.  We refer to this 
mixing as an action of $B$ on $A$.  The following definition is due to 
Kontsevich \cite{kontsevich-operads-motives}.
\end{example}

\begin{definition}
\label{definition:ed-algebra-action-on-edone-algebra}
 Let $B$ be an $E_d$ algebra and $A$ an $E_{d-1}$ algebra.  A \emph{swiss 
 cheese action of $B$ on $A$} is the structure of a swiss cheese algebra on the 
 pair $(B,A)$ extending the given $E_d$ and $\edone$ structures.  We may also 
 simply call this ``an action of $B$ on $A$''.
\end{definition}
\subsection{Statement of main theorem}
Now we can informally state the conjecture proven in this paper.

\begin{theorem*}[Informal statement of Swiss Cheese Conjecture]
The Hochschild cochain object $\hoch(A)$ of an $\edone$ algebra $A$ is the 
universal $E_d$ algebra acting on $A$.  In other words, for any $\edone$ 
algebra $A$, there is an $E_d$ algebra structure on $\hoch(A)$ such that for 
any $E_d$ algebra $B$, giving a map of $E_d$ algebras $B \to \hoch(A)$ is 
equivalent to giving the structure of an $\swcheese$ algebra on the pair 
$(B,A)$ extending the given $E_d$ and $\edone$ structures.
\end{theorem*}

The basic structure of the proof is outlined in 
diagram~\eqref{diagram:summary-of-proof}.  The categories of operads shown in the 
diagram are defined precisely in 
section~\ref{section:other-viewpoints-on-swiss-cheese}. The category 
$\topl^\Sigma$ consists of symmetric sequences of topological spaces.  
Informally, $\Operads(\Coll{\topl})$ consists of $K$-colored operads whose 
$\e$-colored output is trivial, and $\Operads(\Collleqone{\topl})$ further 
restricts to those pieces whose $\e$-colored inputs total 0 or 1.
The forgetful functors are presented below as straight arrows; there is a left 
adjoint shown as a bent arrow.

Boardman and Vogt's $W$ construction \cite{boardman-vogt-73}  is an explicit 
cofibrant replacement functor, which we apply to $\swcheese$, to get an 
equivalent cofibrant operad $\wswcheese\coloneqq W\swcheese$.   The $\W$ 
construction does not strictly commute with the forgetful functors in 
\eqref{diagram:summary-of-proof}. In particular $\W(\swcheese^\h)$ is not 
isomorphic to $(\W \swcheese)^\h$, but they are homotopy equivalent.  We will 
let $\sch$ denote $\W(\swcheese^\h)$ and $\schone = W(\swcheese^{\h 1})$.
\begin{equation}
\label{diagram:summary-of-proof}
\begin{tikzpicture}[numbered picture]
  \matrix (m) [matrix of math nodes, row sep={1.25ex}, column sep={8ex}, text 
  height=1.5ex, text depth=0.25ex]
{ \Op{K}{} & \Op{}{}(\Coll{\topl}) & \Operads(\Collleqone{\topl}) \\
  \wswcheese & \sch & \schone \\
  \schinf\rtimes \E & \schinf & \\};
\path[->]
 (m-1-1) edge (m-1-2)
 (m-1-2) edge (m-1-3)
 (m-1-3) edge[bend right] (m-1-2)
 (m-2-1) edge[|->] (m-2-2)
 (m-2-2) edge[|->] (m-2-3)
 (m-2-3) edge[|->] (m-3-2)
 (m-3-2) edge[dashed] node[above, map name] {transfer of} node[below, map name] 
 {structure} (m-3-1);
 \end{tikzpicture}
\end{equation}
Proposition~\ref{proposition:hoch-A-in-terms-of-swiss-cheese} shows that we can 
use $\schone$ to construct a model for the Hochschild cochain object.  This 
allows us to prove a weak version of the swiss cheese theorem in 
proposition~\ref{proposition:degree0-1-swiss-cheese-theorem} taking place in 
the context of $\Operads(\Collleqone{\topl})$.  Next, we will take the free 
extension of $\schone\in \Operads(\Collleqone{\topl})$ to an operad in 
$\Operads(\Coll{\topl})$, to get $\schinf$.  This immediately gives a version 
of the swiss cheese theorem in the context of $\Operads(\Coll{\topl})$, see 
corollary~\ref{corollary:AlgAschinf(C)-iso-C/Hoch(A)}.  Then we use the fact 
that $\schinf$ is freely generated by its degree 0 and 1 pieces to prove a 
version of the swiss cheese theorem in the context of $\Op{K}{}$.  None of 
these three versions of the swiss cheese theorem make any use of $E_d$.  One 
can think of this last ``universal cheese'' theorem 
(proposition~\ref{proposition:schinf-rtimes-O}) as a construction of the 
universal $K$-colored operad built from $\wedone$ and controlling 
$\wedone$-linear actions on $A$.

To bring $E_d$ back in to the story, we use a technical result, 
\ref{theorem:acyclic-cofibration}, which shows that the canonical map $\schinf 
\to \sch$ is an equivalence.  Observe that one can view $\wswcheese$ as $\sch$ 
equipped with the extra structure of a right action of $E_d$.  Now use a 
transfer of structure argument to construct an operad $\E$ which is equivalent 
to $E_d$ and which acts on the right on $\schinf$.  This allows us to define 
$\schinf \rtimes \E$, which we show is equivalent to $\swcheese$ in 
\ref{lemma:schinf-rtimes-E-simeq-swcheese}.  The universal property of $\schinf 
\rtimes\E$ with respect to the Hochschild cochain object is stated in 
theorem~\ref{theorem:precise-statement-of-main-theorem} and follows from 
proposition~\ref{proposition:schinf-rtimes-O}.

\begin{theorem}[Precise version of Kontsevich's Swiss Cheese 
	Conjecture]\label{theorem:precise-statement-of-main-theorem}
Let $\calC$ be a symmetric monoidal model category tensored over $\topl$ 
and satisfying the conditions in notation~\ref{notation:mathcal-C}. Let $A \in 
\Alg_{\wedone}(\calC)$ be cofibrant and fibrant using the projective model 
structure (\ref{definition:Quillen-adjunctions}).  There is a model of 
$\swcheese$, called $\schinf\rtimes\E$ where $\E\simeq E_d$.  There is also a 
model for the Hochschild cochain object of $A$, called $\hoch(A)$, such that 
$\hoch(A)$ is the universal $\E$ algebra acting on $A$ through $\schinf\rtimes 
\E$.  That is, $\hoch(A)$ is an $\E$ algebra and this structure, together with 
the $\edone$ algebra structure on $A$, can be extended to an $\schinf \rtimes 
\E$ algebra structure on $(\hoch(A), A)$ in such a way that there is an 
isomorphism of categories
\[
	\Alg_{(\schinf \rtimes \E)}^A(\calC) \cong \Alg_{\E}(\mathcal 
	C)_{/\hoch(A)}.
\]
The category on the left consists of $\E$ algebras $B$ together with an action 
of $B$ on $A$.  The category on the right consists of $\E$ algebras $B$ 
together with an $\E$ algebra map $B \to \hoch(A)$.
\end{theorem}

Lemma~\ref{lemma:schinf-rtimes-E-simeq-swcheese} shows that $\schinf \rtimes 
\E$ is equivalent to $\swcheese$. Lemma~\ref{lemma:E-simeq-Ed} shows that $\E$ 
is equivalent to $E_d$.  
Proposition~\ref{proposition:hoch-A-in-terms-of-swiss-cheese} sets up our 
choice of model of $\hoch(A)$.  Proposition~\ref{proposition:schinf-rtimes-O}, 
together with the construction of the operad $\schinf \rtimes \E$ in 
section~\ref{section:right-Ed-action} shows the desired isomorphism of 
categories.

\subsection{Defining \texorpdfstring{$\swcheese^{\h}$}{SCh} and 
\texorpdfstring{$\swcheese^{\h 1}$}{SCh1}}
\label{section:other-viewpoints-on-swiss-cheese}
Recall that $\calC$ is a symmetric monoidal category enriched over 
$\mathcal S$, our basic category in which our operads live.  We will assume 
that both $\calC$ and $\mathcal S$ have all coproducts and that tensor 
products distribute over finite coproducts. In the case of operads from 
$\Op{K}{}$ we have $(\mathcal S, \otimes)=(\topl, \times)$, the symmetric 
monoidal category of compactly generated topological spaces with cartesian 
product.  In the case of operads such as $\swcheese^\h$ and $\sch$ from 
$\Operads(\Coll{\topl})$, we have $(\mathcal S, \otimes) = (\Coll{\topl}, 
\otimes)$ from definition~\ref{definition:symmetric sequences}.  Finally, for 
$\schone$ and $\swcheese^{\h 1}$ we use $\mathcal S=\Collleqone{\topl}$ as in 
definition~\ref{definition:degree-0-1-symmetric sequences}.

\begin{definition}
\label{definition:symmetric sequences}
Let $\Sigma$ denote the \emph{opposite} of the category of finite sets with 
morphisms given by bijections, and let $(\calD, \otimes_\calD)$ be any 
symmetric monoidal category.  The category of functors $\Sigma \to \calD$, 
denoted $\Coll{\calD}$, is usually called the category of \emph{symmetric 
sequences} in $\calD$.  We endow $\Coll{\calD}$ with the usual symmetric 
monoidal structure given by left Kan extension of $\Sigma \times \Sigma \to 
\calD \times \calD \to \calD$ along disjoint union of sets $\Sigma \times 
\Sigma \to \Sigma$.  Specifically, if $X, Y \in \Coll{\calD}$, then $X \otimes 
Y \in \Coll{\calD}$ satisfies
\[
(X \otimes Y) (\underline n) = \coprod_{ n = n_1 + n_2} \ind_{\aut(\underline 
n_1) \times \aut(\underline n_2)}^{\aut(\underline n)} (X(\underline n_1) 
\otimes_{\calD} Y(\underline n_2)),
\]
See Harper \cite[definition 3.3]{harper-2007} for more details.  
\end{definition}

\begin{definition}
\label{definition:degree-0-1-symmetric sequences}
Let $\Sigma^{\leq 1}$ denote the full subcategory of $\Sigma$ consisting of 
finite sets of size 0 or 1 together with bijections as morphisms.  Let 
$\Collleqone{\calD}$ denote the symmetric monoidal category of functors 
$\Sigma^{\leq 1} \to \calD$, with monoidal structure inherited from 
$\Coll{\calD}$.  Call these the \emph{degree 0-1 symmetric sequences in 
$\calD$}. Concretely,  $\Collleqone{\calD}$ is just the category $\calD \times 
\calD$ endowed with the symmetric monoidal structure
\begin{equation}
\label{equation:Coll-leq-1-tensor-product}
(C_0, C_1) \otimes ( D_0, D_1) = (C_0 \otimes_{\calD} D_0, (C_0 \otimes_{\calD} 
D_1) \coprod (C_1 \otimes_{\calD} D_0)).
\end{equation}
The braiding isomorphism $(C_0, C_1) \otimes (D_0, D_1) \to (D_0, D_1) \otimes 
(C_0, C_1)$ is induced from the braiding isomorphism on $\calD$.  
\end{definition}

\begin{definition}
\label{definition:opK}
  Let $\Op{K}{}$ denote the category of $K =\{\e,\h\}$-colored operads in 
  $\topl$.   There is a forgetful functor $\Op{K}{} \to 
  \Operads(\Coll{\topl})$.  This functor takes $\mathcal O \in \Op{K}{}$ to the 
  operad $\calO^\h$ whose arity $m$ component is the symmetric sequence $n 
  \mapsto \calO^\h(n,m)$ (see notation~\ref{notation:underline-notation}).  We 
  forget the spaces $\calO^e(n,m)$, and think of elements of $\calO^h(n,m)$ as 
  degree $n$, arity $m$ elements of $\calO^h$.

  The functor $\Operads(\Coll{\topl}) \to \Operads(\Collleqone{\topl})$ is 
  induced by the symmetric monoidal forgetful functor $\Coll{\topl} \to 
  \Collleqone{\topl}$.  Denote the image of $\calO$ in 
  $\Operads(\Collleqone{\topl})$ as $\calO^{\h 1}$.
\end{definition}

\begin{example}
  \label{example:outline-sch-as-operad-in-coll}
  We outline the structure of $\sch$ as an operad in $\Operads(\Coll{\topl})$.  
  Think of $\sch(\smallbullet,m)$ as the symmetric sequence $n \mapsto 
  \sch(n,m)$.  In \eqref{equation:sch(-,m)-sch(-,ell)}, $\otimes$ is the tensor 
  product of symmetric sequences.  The operad composition law is
  \begin{equation}\label{equation:sch(-,m)-sch(-,ell)}
	\sch(\smallbullet,m) \otimes \sch(\smallbullet, \ell_1) \otimes \cdots 
	\otimes \sch(\smallbullet, \ell_m) \to \sch(\smallbullet, \ell),
  \end{equation}
  where $\ell = \sum_i \ell_i$.  The degree $n$ component of the right hand 
  symmetric sequence is $\sch(n,\ell)$.  The degree $n$ component of the left 
  hand symmetric sequence is
  \[
	\coprod_{n_0 + \cdots + n_{m} = n} \ind\,\, \sch(n_0,m) \times 
	\sch(n_1,\ell_1) \times \cdots \times \sch(n_m, \ell_m),
  \]
  where $\ind$ is the induction functor giving the correct symmetric group 
  action.  The point is that if we delete all appearances of $E_d$ from 
  \eqref{equation:swiss-cheese-full-composition} then it provides exactly the 
  data of \eqref{equation:sch(-,m)-sch(-,ell)}.
\end{example}

\begin{example}
\label{example:schone-as-1-colored-operad}
 We outline the structure of $\swcheese^{\h 1}$ as a 1-colored operad in 
 $\Collleqone{\topl}$. For each $m',m$, define an object $\swcheese^{\h 
 1}(m';m) \in \Collleqone{\topl}$,
\begin{align}
\swcheese^{\h 1}(m';m) \nonumber
&=
(\swcheese^\h(0,m'; 0,m),\swcheese^\h(1, m'; 0, m))\\
\label{equation:edone-identification}
&\cong
(E_{d-1}(m';m),\swcheese^\h(1, m';0,m)) .
\end{align}
  Using the symmetric monoidal structure from 
  definition~\ref{definition:degree-0-1-symmetric sequences} and the 
  identification in equation~\eqref{equation:edone-identification} we can write 
  the operad structure maps on $\swcheese^{\h 1}$ as a triple of morphisms.  
  The map of degree  0 pieces,
 \[
 E_{d-1}(m''; m) \times E_{d-1}(m';m'') \to E_{d-1} (m';m),
 \]
 and the maps of degree 1 pieces,
 \[
 E_{d-1}(m''; m) \times \swcheese^{\h}(1,m';0,m'') \to 
 \swcheese^{\h}(1,m';0,m),
 \]
 and
 \[
 \swcheese^\h(1,m''; 0,m) \times E_{d-1} (m'; m'') \to \swcheese^\h(1,m';0,m).
 \]
\end{example}

\section{Hochschild cohomology from swiss cheese}
\label{section:hochschild-cohomology-from-swiss-cheese}

For the remainder of the paper we replace $\edone, E_d,$ and $\swcheese$ by 
cofibrant models given by the Boardman-Vogt $W$ construction 
\cite{boardman-vogt-73}. A proof that this gives a cofibrant replacement for 
certain operads can be found in Berger-Moerdijk \cite{berger-moerdijk-w}.  We 
will denote these cofibrant replacements by $\wedone, \wed,$ and $\wswcheese$.
We also want to restrict our attention to swiss cheese algebras in categories 
where we can do homotopy theory.  In the proper context the Hochschild cochain 
object of an $\wedone$ algebra $A$ has a natural model constructed from $A$ and 
the degree 0-1 parts of $\wswcheese^\h$. This is the content of 
proposition~\ref{proposition:hoch-A-in-terms-of-swiss-cheese}.
\subsection{Homotopy theoretic context}
\label{section:homotopy-theoretic-context}
\begin{definition}
\label{definition:monoidal-model-category}
From \cite[definition 4.2.6]{hovey-book}, a \emph{symmetric monoidal model 
category} $\mathcal S$ is a closed symmetric monoidal category whose monoidal 
structure $\otimes\colon \mathcal S \times \mathcal S \to \mathcal S$ is a 
Quillen bifunctor, and where the cofibrant replacement $Q1 \to 1$ of the 
monoidal unit induces weak equivalences $Q1 \otimes X \to X$ for every 
cofibrant $X$.
\end{definition}
\begin{example}
\label{example:topl-as-symmetric-monoidal-model-category}
The category $(\topl, \times)$ of compactly generated spaces with the cartesian 
product and Serre model structure is a symmetric monoidal model category.
\end{example}
\begin{definition}
\label{definition:category-over-S}
Let $\mathcal S$ be a symmetric monoidal model category.
A \emph{symmetric monoidal model category tensored over $\mathcal S$} is a 
closed symmetric monoidal model category $\calC$, together with a 
symmetric monoidal Quillen functor $\mathcal S \to \calC$.  For more 
details see Hovey \cite[definition 4.2.20]{hovey-book}.

In particular, $\calC$ comes equipped with functors
\begin{gather*}
\otimes: \calC\times \calC \rightarrow \calC  \quad  
\otimes:\mathcal S \times \calC \rightarrow \calC. \\
\inthom\subc: \calC^{op}\times \calC \rightarrow \mathcal 
C\quad\map\subc: \calC^{op}\times \calC \rightarrow\mathcal S
\end{gather*}
The mapping spaces $\map\subc(A,B)$ give $\calC$ the structure of a 
category enriched over $\mathcal S$, so we can speak of $\wedone, \wed$ and 
$\wswcheese$ algebras in $\calC$.

For any object $A$ of $\calC$, the functor $- \otimes A$ has right 
adjoints $\inthom\subc( A, - ): \calC \rightarrow \calC$ and 
$\map\subc( A, - ): \calC \rightarrow\mathcal S$.
This data satisfies Quillen's SM7 axiom (Hovey \cite[section 4.2]{hovey-book}).  
\end{definition}

\begin{example}
\label{example:Coll-leq-1-tensored-over-Coll-leq-1-topl}
The category $\Collleqone{\calC}$ from 
definition~\ref{definition:degree-0-1-symmetric sequences} is tensored over the 
symmetric monoidal category $\Collleqone{\mathcal \topl}$ with
\[
	\Collleqone{\topl} \otimes \Collleqone{\calC} \to \Collleqone{\mathcal 
	C}
\]
given by the analogue of equation~\eqref{equation:Coll-leq-1-tensor-product},
\[
(X_0, X_1) \otimes (C_0, C_1) = (X_0 \otimes C_0, X_0 \otimes C_1 \coprod X_1 
\otimes C_0),
\]
where $X_i \in\topl$ and $C_i \in \calC$.
\end{example}
\begin{example}
\label{example:Coll-leq-1-tensored-over-topl}
The symmetric monoidal functor $\mathcal \topl \to \Collleqone{\topl}$ sending 
$X$ to $(X, \emptyset)$ makes both $\Collleqone{\topl}$ and 
$\Collleqone{\calC}$ into symmetric monoidal model categories tensored 
over $\topl$.
Note that for a topological 1-colored operad $\mathcal O$, we can consider 
algebras over $\mathcal O$ in $\calC$ as well as algebras over $\mathcal 
O$ in $\Collleqone{\calC}$.
\end{example}
The category of degree 0-1 symmetric sequences is naturally home to $\mathcal 
O$-algebras $A$ and $\mathcal O \dash A$ modules $M$.
\begin{definition}
\label{definition:O-A-modules}
Suppose $\mathcal O$ is a 1-colored operad in $\topl$, and let $(A,M)$ be an 
object of $ \Collleqone{\calC}$.  The structure of an $\calO$ algebra on 
the degree 0-1 symmetric sequence $(A,M)$ is the structure of an $\mathcal O$ 
algebra on $A$ together with the data of maps
\begin{equation*}
\mathcal O (m) \otimes M \otimes A^{\otimes m-1} \to M,
  \end{equation*}
  satisfying certain conditions (see diagram~\eqref{diagram:edone-A-module}).  We 
  call this data the structure of an \emph{$\mathcal O\dash A$ module} on $M$.
 Given a fixed $\mathcal O$ algebra $A \in \calC$, the category 
 $\modules_{\mathcal O}^A(\calC)$ of $\mathcal O\dash A$ modules has 
 objects $M \in \calC$ together with the structure of an $\mathcal O$ 
 algebra on $(A,M)\in\Collleqone{\calC}$ extending the given $\mathcal O$ 
 algebra structure on $A$.

We can enrich $\modules_{\mathcal O}^A(\calC)$ over $\calC$ and over 
$\topl$.
Indeed, given $\mathcal O\dash A$ modules $M'$ and $M$, we can define the 
hom-object of $\mathcal O\dash A$ module morphisms from $M'$ to $M$ as the 
equalizer
\[
\inthom_{{\mathcal O}}^A(M', M) \xrightarrow{\mbox{\scriptsize eq}} 
\inthom_{\calC}(M', M) \rightrightarrows \inthom_{\calC}(F_{\mathcal 
O}^A (M'), M).
 \]
Here $F_{\mathcal O}^A \colon \calC \rightarrow \calC$ is the free $\mathcal 
O\dash A$ module monad on $\calC$.  One can define this using the monad 
$F_\calO$ on $\calC \times \calC$ which sends $(A,M)$ to $(F_\calO^0(A), 
F_\calO^1(A,M))$, where $F_\calO^0$ is the free $\calO$-algebra monad and
\[
	F_\calO^1(A,M) = \coprod_{m} \calO(m)\otimes_{S_m} \Big( \coprod_{i=1}^m 
	A^{\otimes i-1} \otimes M \otimes A^{\otimes m-i} \Big).
\]
The natural transformation $F_\calO^1(F_\calO^0(-), F_\calO^1(-,-)) \to 
F_\calO^1(-,-)$ is given by composition in $\calO$.  This, along with 
$F_\calO^0 F_\calO^0\to F_\calO^0$, defines the monad structure $F_\calO 
F_\calO \to F_\calO$.  Finally, we define the free $\calO\dash A$ algebra monad 
via the coequalizer
\[
	F_\calO^1(F_\calO^0(A), F_\calO^1(A,M)) \rightrightarrows F_\calO^1(A,M) 
	\xrightarrow{\mbox{\scriptsize coeq}} F_\calO^A(M).
\]
 It is clear that every $\mathcal O\dash A$ module $M$ is equipped with a 
 canonical map $F_{\mathcal O}^A (M) \to M$.  The two parallel arrows in the 
 equalizer are given by the two maps $F_{\mathcal O}^A( M') \to M'$ and 
 $F_{\mathcal O}^A( M) \to M$.  The topological space $\map_{\mathcal 
 O}^{A}(M',M)$ is defined as an equalizer in exactly the same manner.
\end{definition}

\begin{example}
\label{example:A-as-wedone-A-module}
Let $A \in \calC$ be an $\wedone$ algebra, then the degree 0-1 symmetric 
sequence $(A,A)$ is naturally an $\wedone$ algebra.  That is, $A$ is naturally 
an $\wedone\dash A$ module.
\end{example}


\begin{definition}
\label{definition:Quillen-adjunctions}
Recall from example~\ref{example:schone-as-1-colored-operad} that 
$\swcheese^{\h 1}$ is an operad in $\Collleqone{\topl}$.  Let $\schone \in 
\Collleqone{\topl}$ be $W(\swcheese^{\h 1})$.  Since $\Collleqone{\calC}$ 
is enriched over $\Collleqone{\topl}$ by 
example~\ref{example:Coll-leq-1-tensored-over-Coll-leq-1-topl}, we can consider 
$\schone$ algebras in $\Collleqone{\calC}$.  In addition, by  
example~\ref{example:Coll-leq-1-tensored-over-topl} we can consider $\wedone$ 
algebras in $\Collleqone{\calC}$.
There are adjunctions
\begin{gather*}
\calC
\leftrightarrows
\Alg_{\wedone}(\calC)
\quad
\calC
\leftrightarrows
\modules_{\wedone}^A(\calC)
\\
\calC \times \calC
\leftrightarrows
\Alg_{\wedone}(\Collleqone{\calC})
\leftrightarrows
\Alg_{\schone}(\Collleqone{\calC})
\end{gather*}
We describe the right adjoints only.  On the top left the $\wedone$ algebra $A$ 
is sent to the underlying object $A$ of $\calC$.  The top right functor 
sends the $\wedone\dash A$ module $M$ to the underlying object $M$ of $\mathcal 
C$.  In the pair of composable  adjunctions, an $\schone$ algebra $(A,M)$ can 
be considered as an $\wedone$ algebra by forgetting the structure maps in 
equation~\eqref{equation:schone-alg-structure-0-1}.  For the final adjunction, 
any $\wedone$ algebra in $\Collleqone{\calC}$ has an underlying pair of 
$\calC$-objects $(A,M)$.
\end{definition}

\begin{notation}
\label{notation:mathcal-C}
Throughout the remainder of this paper $\calC$ will be a cofibrantly 
generated symmetric monoidal model category tensored over $\mathcal \topl$ such 
that the adjunctions in definition~\ref{definition:Quillen-adjunctions} are 
Quillen adjunctions.
\end{notation}

\begin{remark}
	In \cite[theorem 6]{spitzweck-2001}, Spitzweck shows that 
	$\modules_{\mathcal O}^A(\calC) \leftrightarrows \calC$ is a Quillen 
	adjunction if $A$ is cofibrant in $\calC$.  In addition
		Berger-Moerdijk show in \cite[proposition 
		4.1]{berger-moerdijk-axiomatic} that if $\calD$ is a symmetric monoidal 
		model category which is cofibrantly generated, has a cofibrant unit, 
		and a symmetric monoidal fibrant replacement functor, then the category 
		$\Alg_{\calO}(\calD)$ has the projective model structure induced from 
		the forgetful functor to $\calD$.  The operad $\calO$ in this theorem 
		is an operad in $\calD$.  However, their result is more general, as 
		seen in remark 4.6.4 of the same article.  Their argument extends 
		without change to show that $\Alg_{\calO} (\calD)$ has the desired 
		model structure in the case that $\calO$ is an operad in $\topl$ and 
		$\calD$ is tensored over $\topl$ as in definition 
		\ref{definition:monoidal-model-category}, and the generating trivial 
		cofibrations of $\calD$ are cofibrant.  The condition that the monoidal 
		unit of $\calD$ is cofibrant is not necessary in this situation.  The 
		cofibrance of the monoidal unit in $\topl$ is enough.  Taking $\calO = 
		\wedone$ or $\calO = \schone$ and $\calD = \calC$ or $\calD = 
		\Collleqone{\calC}$, we conclude that the assumptions in notation 
		\ref{notation:mathcal-C}, and in the main theorem of this paper, hold 
		when $\calC$ is the category of compactly generated Hausdorff spaces 
		$(\topl)$, or spectra, or chain complexes.
	\end{remark}
We use the model structure on the category of $\wedone\dash A$ modules to 
define the Hochschild cochain object.
\begin{definition}
\label{definition:hochschild-cohomology}
Given an $\wedone$ algebra $A \in \calC$, let the \emph{Hochschild cochain 
object} of $A$ be
\[
 \hoch(A) = \inthom_{\wedone}^{A} (A^{c}, A^{f}),
\]
where $\inthom_{\wedone}^{A}$ is given by the equalizer in 
definition~\ref{definition:O-A-modules}. The $\wedone\dash A$ modules $A^c$ and 
$A^f$ are cofibrant and fibrant replacements for $A$ respectively.
Note that $\hoch(A)$ is an object of $\calC$.
\end{definition}

We will use the degree 0 and 1 pieces of the swiss cheese operad to build a 
cofibrant replacement for $A$ as an $\wedone\dash A$ module.

\subsection{Swiss cheese in degrees zero and one}
\label{section:swiss-cheese-in-degree-0-1}
 An $\schone$ algebra (see definition~\ref{definition:Quillen-adjunctions} and 
 example~\ref{example:schone-as-1-colored-operad}) is a pair $(A,M)$ of objects 
 of $\calC$ together with maps in $\Collleqone{\calC}$ for every $m$,
\begin{equation}
\label{equation:schone-alg-structure-pair-maps}
\big(\wedone(m), \sch(1,m)\big) \otimes (A,M) ^{\otimes m} \to (A,M),
\end{equation}
where we have used the isomorphism $\sch(0,m) \simeq \wedone(m)$.
Alternatively, we can view the morphism in 
equation~\eqref{equation:schone-alg-structure-pair-maps} as three separate maps 
in $\calC$.
\begin{align}
\label{equation:schone-alg-structure-0-0}
&\wedone(m) \otimes A ^{\otimes m} \to A\\
\label{equation:schone-alg-structure-1-1}
& \wedone(m) \otimes M \otimes A^{\otimes m - 1} \to M \\
\label{equation:schone-alg-structure-0-1}
& \sch(1,m) \otimes A^{\otimes m} \to M
\end{align}
The condition that the maps in 
equation~\eqref{equation:schone-alg-structure-pair-maps} define an $\schone$ 
structure on the pair $(A,M)$ is the condition that the 
diagram~\eqref{diagram:schone-alg-pair} commutes in the 
$\Collleqone{\topl}$-enriched category $\Collleqone{\calC}$.
\begin{equation}
\label{diagram:schone-alg-pair}
\begin{tikzpicture}[numbered picture]
  \matrix (m) [matrix of math nodes, row sep={3ex}, column sep={1em}, text 
  height=1.5ex, text depth=0.25ex]
{|(a)| \hspace{-0.2in} 
\schone(m') \otimes \schone(m;m') \otimes (A,M)^{\otimes m}
  &
  |(b)|
  \schone(m) \otimes (A,M)^{\otimes m} \hspace{-.2in}
  \\
|(c)|
\schone(m') \otimes (A,M)^{\otimes m'}
&
|(d)|
(A,M)
\\};
 \path[->]
 (a) edge (b)
 (a) edge (c)
 (b) edge (d)
 (c) edge (d);
 \end{tikzpicture}
\end{equation}
In terms of equations~\eqref{equation:schone-alg-structure-0-0}, 
\eqref{equation:schone-alg-structure-1-1}, and 
\eqref{equation:schone-alg-structure-0-1}, diagram~\eqref{diagram:schone-alg-pair} 
splits into four diagrams.  Each diagram is determined by the degrees of the 
three tensor factors in the upper left hand corner of 
diagram~\eqref{diagram:schone-alg-pair}.
In the first the degrees are $0, 0, 0$; in the second the degrees are $0, 0, 
1$; in the third, $1, 0, 0$; and in the fourth, $0, 1, 0$.
\begin{equation}
\label{diagram:edone-algebra}
\begin{tikzpicture}[numbered picture]
  \matrix (m) [matrix of math nodes, row sep={3ex}, column sep={5ex}, text height=1.5ex, text depth=0.25ex]
{\wedone(m') \otimes \wedone(m;m') \otimes A^{\otimes m}
& \wedone(m) \otimes A^{\otimes m} \\
\wedone(m') \otimes A^{\otimes m'}
& A \\};
\path[->]
 (m-1-1) edge (m-1-2)
 (m-1-1) edge (m-2-1)
 (m-2-1) edge (m-2-2)
 (m-1-2) edge (m-2-2);
 \end{tikzpicture}
\end{equation}
The above diagram,~\eqref{diagram:edone-algebra}, commutes for all $m', m$ if and only if $A$ is an $\wedone$ algebra.
Diagram \eqref{diagram:edone-A-module} below commutes if and only if $M$ is an 
$\wedone\dash A$ module.
\begin{equation}
\label{diagram:edone-A-module}
\begin{tikzpicture}[numbered picture]
  \matrix (m) [matrix of math nodes, row sep={3ex}, column sep={1em}, text height=1.5ex, text depth=0.25ex]
{\hspace{-.1in}
\wedone(m') \otimes \wedone(m;m') \otimes M \otimes A^{\otimes m-1}
&
\wedone(m) \otimes M  \otimes A^{\otimes m'-1}
\hspace{-.1in}
\\
\wedone(m) \otimes M \otimes A^{\otimes m-1}
& M \\};
\path[->]
 (m-1-1) edge (m-1-2)
 (m-1-1) edge (m-2-1)
 (m-2-1) edge (m-2-2)
 (m-1-2) edge (m-2-2);
 \end{tikzpicture}
\end{equation}
Diagram \eqref{diagram:sch-edone-compatibility} shows a compatibility condition between the degree 0 and degree 1 structures.
\begin{equation}
\label{diagram:sch-edone-compatibility}
\begin{tikzpicture}[numbered picture]
  \matrix (m) [matrix of math nodes, row sep={3ex}, column sep={2em}, text height=1.5ex, text depth=0.25ex]
{\sch(1,m') \otimes \wedone(m;m') \otimes A^{\otimes m}
&\sch(1,m) \otimes A^{\otimes m} \\
\sch(1, m') \otimes A^{\otimes m'}
& M \\};
\path[->]
 (m-1-1) edge (m-1-2)
 (m-1-1) edge (m-2-1)
 (m-2-1) edge (m-2-2)
 (m-1-2) edge (m-2-2);
 \end{tikzpicture}
\end{equation}
Diagram \eqref{diagram:edone-sch-compatibility} presents another compatibility condition between the degree 0 and degree 1 structures.
\begin{equation}
\label{diagram:edone-sch-compatibility}
\begin{tikzpicture}[numbered picture]
  \matrix (m) [matrix of math nodes, row sep={3ex}, column sep={2em}, text height=1.5ex, text depth=0.25ex]
{\wedone(m')  \otimes \sch(1, m; 0,m') \otimes A^{\otimes m}
&\sch(1, m) \otimes A^{\otimes m} \\
\wedone(m') \otimes M\otimes A^{\otimes m'-1}
& M \\};
\path[->]
 (m-1-1) edge (m-1-2)
 (m-1-1) edge (m-2-1)
 (m-2-1) edge (m-2-2)
 (m-1-2) edge (m-2-2);
 \end{tikzpicture}
\end{equation}

\begin{figure}[ht!]
  \begin{center}
	\includegraphics{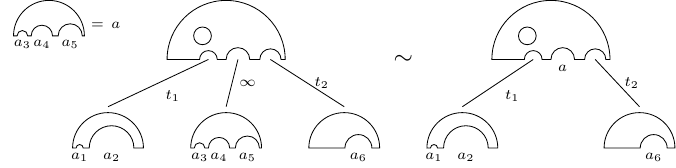}
	\caption[Relations in $A^{sc}$]{The relations in $A^{sc}$ come from the 
	$\wedone$ algebra structure of $A$.  If $m$ is the map $A^{\otimes 3} \to 
	A$ given by the swiss cheese element in $\wswcheese^\h(0,3)\simeq 
	\wedone(3)$ in the figure, set $a = m(a_1, a_2, a_3)$.  The edges $t_1$ 
	and $t_2$ are less than $\infty$, so the relation does not apply to the 
	vertices on the left and right.}
	\label{figure:A-sc-equivalence-relation}
  \end{center}
\end{figure}

\begin{example}
\label{example:A-sc-as-schone-algebra}
In this example we construct the universal extension of an $\wedone$ algebra 
$A$ to an $\schone$ algebra.  We denote this universal pair by $(A, A^{sc})$.
The composite forgetful functor
\[
\Alg_{\schone}(\Collleqone{\calC})
\to
\Alg_{\wedone}(\Collleqone{\calC})
\to
 \Alg_{\wedone}(\calC),
\]
has a left adjoint which sends the $\wedone$ algebra $A$ to the pair $(A, 
A^{sc})$ where $A^{sc}$, which may be read as ``$A$ swiss cheese'', is a 
quotient of
\begin{equation}
\label{equation:bar-A-sc}
 \bar{A}^{sc} = \coprod_{m\geq 0} \sch(1,m) \otimes_{S_m} A^{\otimes m},
\end{equation}
where $S_m = \aut(\ul m)$.
We can think of $\bar{A}^{sc}$ heuristically as $\sch(1,\smallbullet) \otimes 
A^{\otimes \bullet}$.  Now both
$A$  and
$\sch(1, \smallbullet)$
carry an action of
$\wedone$, so we can form the quotient $A^{sc} \coloneqq 
\sch(1,\smallbullet)\otimes_{\wedone} A^{\otimes \smallbullet}$.  More 
precisely, $A^{sc}$ is defined as the coequalizer
\begin{equation}
\label{equation:A-sc}
 \coprod_{m,m'} \sch(1,m) \otimes \wedone(m';m) \otimes A^{\otimes m'} 
 \rightrightarrows \bar{A}^{sc}\xrightarrow{\mbox{\scriptsize coeq}} A^{sc},
\end{equation}
where one of the arrows is given by the operadic composition on swiss cheese
and the other by the $\wedone$ structure on A.
 Figure~\ref{figure:A-sc-equivalence-relation} shows the relation $\sim$ such 
that $A^{sc} = \bar{A}^{sc}\modsim$.
\end{example}

Verifying that $(A, A^{sc})$ is an $\schone$ algebra is a matter of using the 
commuting left and right actions of $\wedone$ on $\schone$.  By this we mean 
the morphism below uses both left and right actions, and can be obtained by 
performing the left action first, then the right, or vice versa,
\[
\wedone(m'') \times \sch(1,m'; 0,m'') \times \wedone( m; m') \to \sch(1,m).
\]
The left action
defines a map
\begin{equation}
\label{equation:wedone-A-on-bar-A-sc}
\wedone(m)\otimes \bar A^{sc} \otimes A^{\otimes m-1} \to \bar A^{sc},
\end{equation}
using the the $\circ_1$ operad composition.
Since the left and right $\wedone$ actions on $\schone$ commute, the arrow in 
~\eqref{equation:wedone-A-on-bar-A-sc} descends to give the data of a 
$\wedone\dash A$ module structure on $A^{sc}$, i.e.  
equation~\eqref{equation:schone-alg-structure-1-1} for $M = A^{sc}$.  Of course, 
the maps from equation~\eqref{equation:schone-alg-structure-0-1} with $M = 
A^{sc}$ are simply given by $\bar A^{sc} \to A^{sc}$.

Now let us observe that the four 
diagrams~\eqref{diagram:edone-algebra}-\eqref{diagram:edone-sch-compatibility} 
commute for $(A, A^{sc})$.  The first diagram, \eqref{diagram:edone-algebra}, is 
trivial since $A$ is an $\wedone$ algebra.  The second diagram, 
\eqref{diagram:edone-A-module}, commutes since the left action of $\wedone$ on 
$\schone$ is indeed an action.  That is, it is compatible with composition in 
$\wedone$.  The third diagram, \eqref{diagram:sch-edone-compatibility}, clearly 
commutes.  Indeed, this diagram is the reason the coequalizer $\bar A^{sc} \to 
A^{sc}$ in equation~\eqref{equation:A-sc} was defined in the first place.  
Finally, the fourth diagram, \eqref{diagram:edone-sch-compatibility}, certainly 
commutes if $M = \bar A^{sc}$.  In this case, note that the bottom map in 
diagram~\eqref{diagram:edone-sch-compatibility} corresponds to 
equation~\eqref{equation:wedone-A-on-bar-A-sc}.  Thus if we pass from $\bar 
A^{sc}$ to $A^{sc}$ this last diagram still commutes since, by definition, the 
$\wedone \dash A$ module structure on $A^{sc}$ is defined using the quotient 
map $\bar A^{sc} \to A^{sc}$ together with 
equation~\eqref{equation:wedone-A-on-bar-A-sc}.

\begin{lemma}
\label{lemma:mod-sc-leq-1-is-mod-edone-under-A-sc}
Fix an $\wedone$ algebra $A$ and consider $A^{sc}$ as an $\wedone\dash A$ 
module, then there is an isomorphism of categories
\[
\modules_{\schone}^A(\calC) \cong \modules_{\wedone}^A(\mathcal 
C)_{A^{sc}/},
\]
where $\modules_{\schone}^A(\calC)$ is the fiber over $A$ of the forgetful 
functor $\Alg_{\schone}(\Collleqone{\calC})\to \Alg_{\wedone}(\calC)$
\end{lemma}
\begin{proof}
Let $(A,M)$ be an $\schone$ algebra extending the existing $\wedone$ algebra 
structure on $A$.  Then $M \in \modules_{\wedone}^A(\calC)$ and the 
structure maps in equation~\eqref{equation:schone-alg-structure-0-1}, when 
combined for all $m$, give a map $\bar A^{sc} \to M$.  This descends to an 
$\wedone\dash A$ module map $A^{sc} \to M$.

On the other hand if $M$ is an $\wedone\dash A$ module,  then $M$ is already 
equipped with the data of equation~\eqref{equation:schone-alg-structure-1-1}.  If 
$A^{sc} \to M$ is a morphism of $\wedone\dash A$ modules, then $M$ is equipped 
with the data of equation~\eqref{equation:schone-alg-structure-0-1}.  The 
diagram~\eqref{diagram:schone-alg-pair} commutes because of relation defining 
$A^{sc}$ in equation~\eqref{equation:A-sc} and because $A^{sc} \to M$ is a 
morphism of $\wedone\dash A$ modules.
\end{proof}

\begin{corollary}
\label{corollary:A-sc-cofibrant}
Let $A$ be a cofibrant $\wedone$ algebra, then $A^{sc}$ is a cofibrant 
$\wedone\dash A$ module.
\end{corollary}
\begin{proof}
The forgetful functor $\Alg_{\schone}(\calC) \to \Alg_{\wedone}(\mathcal 
C)$ preserves fibrations (see definition~\ref{definition:Quillen-adjunctions}).  
Thus the left adjoint of this forgetful functor, applied to the cofibrant 
object $A$, gives a cofibrant $\schone$ algebra $(A, A^{sc})$.  Thus $A^{sc}$ 
is cofibrant in $\modules^A_{\schone}(\calC)$.
Lemma~\ref{lemma:mod-sc-leq-1-is-mod-edone-under-A-sc} shows that the forgetful 
functor $\modules^A_{\schone}(\calC) \to \modules^A_{\wedone}(\calC)$ 
preserves pushouts.  The model structures here are cofibrantly generated, so 
$A^{sc}$ is also cofibrant as an object of $\modules^A_{\wedone}(\calC)$.
\end{proof}

\begin{definition}
\label{definition:schone-to-edone}
Let $p(m):\sch(1, m)\rightarrow \sch(0,m) \simeq \wedone(m)$ be the projection 
which forgets the single full disk. We can make $p = (p(m))_{m\geq 0}$ into a 
morphism of operads in the following way.  For each $m$ consider the degree 0-1 
symmetric sequence $(\wedone)^{\leq 1}(m) :=(\wedone(m), \wedone(m))$.  The 
structure of $\wedone$ as an operad in $\topl$ can be used to make 
$(\wedone)^{\leq 1}$ an operad in $\Collleqone{\topl}$. This makes $(\id, p) 
\colon \schone \to (\wedone)^{\leq 1}$ into a morphism of operads.

If $A$ is an $\wedone$ algebra we can define a morphism in $\calC$,
\[
\bar A^{sc} = \coprod_{m \geq 0} \sch(1,m)\otimes_{S_m} A^{\otimes m} \to 
\coprod_{m \geq 0} \wedone(m) \otimes_{S_m} A^{\otimes m} \to A,
\]
where the first arrow uses $p$ and the second arrow uses the $\wedone$ algebra 
structure on $A$. This map
factors to give a morphism of $\wedone\dash A$ modules $p_A \colon A^{sc} \to 
A$.
\end{definition}
\begin{remark}
\label{remark:topl-structure-on-operads}
By \cite[section 5]{spitzweck-2001} we can conclude that 
$\Op{}{}(\Collleqone{\topl})$ is tensored over $\topl$.  If $\mathcal O = 
(\mathcal O^0, \mathcal O^1)$ is an operad in degree 0-1 symmetric sequences 
of topological spaces, and $K \in \topl$, then $K \otimes \mathcal O$ is 
defined to be the coequalizer
\[
F(F(K \otimes \mathcal O)) \rightrightarrows F(K \otimes \mathcal 
O)\xrightarrow{\mbox{\scriptsize coeq}} K \otimes \mathcal O,
\]
where $K \otimes \mathcal O$, is the symmetric sequence of degree 0-1 
symmetric sequences whose arity $m$, degree $i$ component is $K \times 
\mathcal O^i(m)$, and $F$ is the free operad functor.

If $\mathcal O$ is a cofibrant operad in $\Collleqone{\topl}$, then $[0,1] 
\otimes \mathcal O$ is a cylinder object, and a homotopy $h \colon [0,1] 
\otimes \mathcal O \to \mathcal P$ gives the data of maps $h^i(m) \colon [0,1] 
\times \mathcal O^i (m) \to \mathcal P^i(m)$, which assemble into $h^i(m'; m) 
\colon [0,1] \times \mathcal O^i(m';m) \to \mathcal P(m';m)$ for $i = 0,1$ and 
$m', m \geq 0$.  The homotopy $h$ is compatible with operad composition in the 
sense that if $\alpha \in \mathcal O^i(m)$ and $\beta \in \mathcal O^j(m'; 
m)$,  $i+j \leq 1$, and $t \in [0,1]$, then
\[
h^i(m)(t, \alpha)\circ h^j( m'; m)(t, \beta) = h^{i+j}(m')(t, \alpha \circ 
\beta).
\]
\end{remark}

\begin{lemma}
\label{lemma:A-sc-equivalent-to-A}
For any $\wedone$ algebra $A \in \calC$, the map $A^{sc} \to A$ is a weak 
equivalence of $\wedone\dash A$ modules.
\end{lemma}

\begin{proof}
Abusing notation we write $p := (\id, p)$ for the morphism of operads from 
definition~\ref{definition:schone-to-edone}.  One can show that $p$ is a weak 
equivalence of fibrant and cofibrant operads in degree 0-1 symmetric sequences 
of topological spaces.  Therefore there is a map of operads $\iota \colon 
(\wedone)^{\leq 1} \to \schone$ and there are homotopies
 $h \colon \id_{\schone} \simeq \iota p$, and $g \colon \id_{(\wedone)^{\leq 
 1}} \simeq p\iota$.
Since $[0,1]\otimes -$ distributes over coequalizers, we can use $h$ to define 
a homotopy $h_A \colon [0,1] \otimes A^{sc} \to A^{sc}$,
\[
\begin{tikzpicture}
  \matrix (m) [matrix of math nodes, row sep={6ex}, column sep={5ex}, text 
  height=1.5ex, text depth=0.25ex]
{\displaystyle \coprod_{m} [0, 1] \otimes\sch(1, m)\otimes A^{\otimes m}
&
\left[0, 1\right] \otimes A^{sc}
\\
\displaystyle \coprod_{m} \sch(1, m)\otimes A^{\otimes m}
&
A^{sc},
\\};
\path[->]
 (m-1-1) edge (m-1-2)
 (m-1-1) edge node[right, map name] {$\sqcup h^1(m) \otimes \id_A^{\otimes 
 m}$} (m-2-1)
 (m-2-1) edge (m-2-2)
 (m-1-2) edge node[right, map name] {$h_A$} (m-2-2);
 \end{tikzpicture}
\]
where $h^1(m)$ is defined from $h$ as in 
remark~\ref{remark:topl-structure-on-operads}.  When $t=1$, $h_A(1, -)$ 
factors as
\[
\begin{tikzpicture}
  \matrix (m) [matrix of math nodes, row sep={6ex}, column sep={5ex}, text 
  height=1.5ex, text depth=0.25ex]
{
\displaystyle \coprod_{m} \sch(1, m)\otimes A^{\otimes m}
&
A^{sc}
\\
\displaystyle \coprod_{m} \wedone(m)\otimes A^{\otimes m}
&
A
\\
\displaystyle \coprod_{m} \sch(1,m) \otimes A^{\otimes m}
&
A^{sc},
\\};
\path[->]
 (m-1-1) edge (m-1-2) edge (m-2-1)
 (m-1-2) edge node[right, map name] {$p_A$} (m-2-2)
 (m-2-1) edge (m-2-2) edge (m-3-1)
 (m-2-2) edge[bend right, dashed] node[above, map name] {$s$} (m-2-1) edge 
 node[right, map name] {$\iota_A$} (m-3-2)
 (m-3-1) edge (m-3-2)
 (m-2-2) edge (m-3-2);
 \end{tikzpicture}
\]
where $p_A$ is the map from definition~\ref{definition:schone-to-edone}, and 
$\iota_A $ is the evident composite in the diagram using $s$.
The map $s$ is a section of the middle horizontal arrow, defined using the 
identity of the operad $\wedone$,
\[
A \simeq \{1_{\wedone}\} \otimes A  \to  \wedone(1) \otimes A \to \coprod_{m}  
\wedone(m) \otimes A^{\otimes m}.
\]
The map $\iota_A$ is a map of $\wedone\dash A$ modules, and $h_A$ is a 
homotopy $\id_{A^{sc}} \simeq \iota_A p_A$.  Similarly, $g$ defines a homotopy 
$g_A \colon \id_{A} \to p_A \iota_A$.
\end{proof}

The precise sense in which Hochschild cohomology can be obtained from the 
degree 0-1 pieces of the swiss cheese operad is contained in the following 
proposition.

\begin{proposition}
\label{proposition:hoch-A-in-terms-of-swiss-cheese}
Let $A$ be a fibrant and cofibrant $\wedone$ algebra.  Then the Hochschild 
cochain object of $A$ can be computed as
\[
\hoch(A) \simeq
 \inthom_{\wedone}^A(A^{sc}, A).
\]
\end{proposition}
\begin{proof}
We are using the projective model structure from 
definition~\ref{definition:Quillen-adjunctions}, so $A$ is fibrant as an object 
of $\calC$ and thus as a $\wedone\dash A$ module.  By 
corollary~\ref{corollary:A-sc-cofibrant} and 
lemma~\ref{lemma:A-sc-equivalent-to-A}, $A^{sc}$ is a cofibrant replacement for 
$A$ as an $\wedone\dash A$ module.  By 
definition~\ref{definition:hochschild-cohomology}, this proves the proposition.
\end{proof}

\section{The universal cheese theorem}
\label{section:universal-cheese}

In the one-colored operad $\schone$ the single full disk was never considered 
as input, only as a marker of degree one.  Allowing the single disk to be 
considered as giving an input means viewing $\schone$ as a partially defined 
$2$-colored operad.  Rather than making the notion \emph{partially defined} 
precise, we simply set up the notion of a 2-colored algebra over $\schone$ in 
definition~\ref{definition:2-colored-schone-algebra}.
Proposition~\ref{proposition:degree0-1-swiss-cheese-theorem} is a version of 
the swiss cheese theorem for the operad $\schone$.  That is, a 2-colored 
$\schone$ structure on the pair $( B, A )$ is equivalent to a $\mathcal 
C$-morphism, $B \rightarrow\hoch ( A )$.  In other words, $\hoch(A)$ is the 
universal object of $\calC$ acting on the $\wedone$ algebra $A$ through 
$\schone$. This result is generalized twice, first in 
corollary~\ref{corollary:AlgAschinf(C)-iso-C/Hoch(A)}, then in 
proposition~\ref{proposition:schinf-rtimes-O}.  We refer to 
proposition~\ref{proposition:schinf-rtimes-O} as the universal cheese theorem 
since it replaces $E_d$ in the swiss cheese theorem with an arbitrary operad.
\addtocounter{subsection}{1}

\begin{definition}\label{definition:2-colored-schone-algebra}
Let $A$ be an $\wedone$ algebra.
Let $\calC_{/\hoch(A)}$ denote the over category of $\hoch(A) \in \mathcal 
C$.
More precisely, the objects are $\calC$-morphisms
$ B \to \hoch(A)$ and the morphisms are $\calC$ morphisms $B \to B'$ 
commuting with the maps to $\hoch(A)$. In addition, let $\Alg_{\schone}^A( 
\calC )$ denote the category of $\schone$ algebras of the form $(B,A)$ 
where the induced $\wedone$ structure on $A$ is the one given.  Morphisms are 
maps of $\schone$ algebras which are identity on $A$.

Given a pair of objects $(B,A)$ of $\calC$, we let $\enD^{\h 1}(B,A)$ denote 
the operad obtained by applying the forgetful functor $\Op{K}{} \to  
\Operads(\Collleqone{\topl})$ to $\enD(B,A)$ from 
\eqref{equation:endomorphism-operad-h} and \eqref{equation:endomorphism-operad-e}.

Let $\Alg^A_{\schone}(\calC)$ denote the category of objects $B \in \calC$ 
together with a morphism in $\Operads(\Collleqone{\topl})$, $\schone \to 
\enD^{\h 1}(B,A)$ extending the $\wedone$ structure on $A$.  A morphism $B \to 
B'$ in $\calC$ induces a map of operads $\enD^{\h 1}(B',A) \to \enD^{\h 
1}(B,A)$.  Such a morphism gives a map in $\Alg^A_{\schone}(\calC)$ if this 
induced map respects to maps from $\schone$.
\end{definition}

\begin{proposition}
\label{proposition:degree0-1-swiss-cheese-theorem}
There is an isomorphism of categories
\[
\Alg^A_{\schone}(\calC) \cong \calC_{/\hoch(A)}.
\]
\end{proposition}
\begin{proof}
The data of an algebra on the left hand side is an object $B \in \calC$ 
together with maps
 \begin{align}
\label{equation:schone-alg-structure-SBA-A}
 \sch(1,m)\otimes B \otimes A^{\otimes m } \to A,
\end{align}
for each $m \geq 0$.
The conditions on \eqref{equation:schone-alg-structure-SBA-A} are that 
diagrams~\eqref{diagram:ESBA-A} and \eqref{diagram:SEBA-A} commute.
\begin{equation}
\label{diagram:ESBA-A}
  \begin{tikzpicture}[numbered picture]
  \matrix (m) [matrix of math nodes, row sep={3ex}, column sep={1em}, text 
  height = 1.5ex, text depth = 0.25ex]
 {
 \wedone(m') \otimes \sch(1,m; 0 , m')\otimes B \otimes A^{\otimes m}
 &
 \wedone( m') \otimes A^{\otimes m'} 
 \\
\sch(1, m) \otimes B\otimes A^{\otimes m}
&
 A,
 \\};
\path[->]
 (m-1-1) edge (m-1-2) edge (m-2-1)
 (m-1-2) edge (m-2-2)
 (m-2-1) edge (m-2-2);
 \end{tikzpicture}
\end{equation}
\begin{equation}
\label{diagram:SEBA-A}
  \begin{tikzpicture}[numbered picture]
  \matrix (m) [matrix of math nodes, row sep={3ex}, column sep={1em}, text 
  height = 1.5ex, text depth = 0.25ex]
 { 
 \sch(1, m')\otimes\wedone(m; m') \otimes B \otimes A^{\otimes m}
 &
 \sch(1, m) \otimes B \otimes A^{\otimes m}
 \\
 \sch(1, m') \otimes B \otimes A^{\otimes m'} 
&
 A.
 \\};
\path[->]
 (m-1-1) edge (m-1-2) edge (m-2-1)
 (m-1-2) edge (m-2-2)
 (m-2-1) edge (m-2-2);
 \end{tikzpicture}
\end{equation}

Equivalently, we can use the hom-tensor adjunction and assemble the maps in 
\eqref{equation:schone-alg-structure-SBA-A} to a single map
$B \to \inthom\subc(\bar A^{sc} , A)$ (see \eqref{equation:bar-A-sc}).  The 
commutativity of diagram~\eqref{diagram:SEBA-A} is equivalent to this map lifting 
to
\begin{equation}
\label{equation:B-homAsc/A}
B \to \inthom\subc( A^{sc} , A).
\end{equation}
Note that, dual to \eqref{equation:A-sc}, $\inthom\subc(A^{sc},A)$ is given by 
the equalizer
\[
\inthom\subc( A^{sc} \!\!, A)
\xrightarrow{\mbox{\scriptsize eq}} \inthom\subc(\bar A^{sc}\!\!,A)
\rightrightarrows
\inthom\subc(\!\coprod_{m,m'} \sch(1,m') \otimes \wedone(m;m') \otimes 
A^{\otimes m}\!\!, A).
\]
With this observation we can now rewrite diagram~\eqref{diagram:ESBA-A} as
\begin{equation}
\label{diagram:BFAsc-A}
  \begin{tikzpicture}[baseline=(current bounding box.center)]
  \matrix (m) [matrix of math nodes, row sep={2ex}, column sep={1em}, text 
  height = 1.5ex, text depth = 0.55ex]
 {
 B \otimes F^A_{\wedone}(A^{sc})
 &
 F^A_{\wedone}(A)
 \\
 B\otimes A^{sc}
&
A.
 \\};
\path[->]
 (m-1-1) edge (m-1-2) edge (m-2-1)
 (m-1-2) edge (m-2-2)
 (m-2-1) edge (m-2-2);
 \end{tikzpicture}
 \end{equation}
Recall that $F_{\wedone}^A\colon \calC \to \calC$ is the free $\wedone\dash A$ 
module functor from definition~\ref{definition:O-A-modules}.  Clearly, 
diagram~\eqref{diagram:ESBA-A} commutes if and only if 
diagram~\eqref{diagram:BFAsc-A} commutes and if and only if the map 
\eqref{equation:B-homAsc/A} factors through $\hoch(A) = 
\inthom_{\wedone}^A(A^{sc},A)$.  We conclude that the data of an $\schone$ 
algebra structure on $(B,A)$ is the data of an $\wedone$ algebra structure on 
$A$ together with a $\calC$-morphism $B \to \hoch(A)$.

It is clear that a map $(B,A) \to (B',A)$ which is identity on $A$ gives an 
$\schone$ algebra morphism if and only if the map $B \to B'$ commutes with the 
corresponding morphisms to $\hoch(A)$.
\end{proof}
\begin{definition}
	\label{definition:schinf}
  Define $\schinf \in \Operads(\Coll{\topl})$ as the image of the left adjoint 
  of $\Operads(\Coll{\topl}) \to \Operads(\Collleqone{\topl})$ applied to 
  $\schone$.  Following definition~\ref{definition:2-colored-schone-algebra}, 
  define $\Alg_{\schinf}^A(\calC)$ as the category of objects $B \in \calC$ 
  together with a morphism in $\Operads(\Coll{\topl})$, $\schinf \to 
  \enD^\h(B,A)$ extending the $\wedone$ structure on $A$.
\end{definition}

\begin{corollary}
  \label{corollary:AlgAschinf(C)-iso-C/Hoch(A)}
  There is an isomorphism of categories
  \[
	\Alg_{\schinf}^A(\calC) \cong \calC_{/\hoch(A)}.
  \]
\end{corollary}
\begin{proof}
	The adjunction isomorphism puts operad maps $\schone \to \enD^{\h1} (B,A)$ 
	in one-to-one correspondence with operad maps $\schinf \to \enD^\h(B,A)$. 
	This gives an isomorphism between $\Alg_{\schinf}^A(\calC)$ and 
	$\Alg_{\schone}^A(\calC)$.  Now apply proposition 
	\ref{proposition:degree0-1-swiss-cheese-theorem}.
\end{proof}

\begin{definition}
  For each $n \geq 0$, let $\schinf(n,\smallbullet)$ denote the operad in 
  $\Operads(\Collleqone{\topl})$ whose arity $m$ component is the degree 0-1 
  symmetric sequence $(\wedone(m), \schinf(n,m))$.  One may think of 
  $\schinf(n, \smallbullet)$ as a bimodule over $\wedone$.
  Let $\enD(\schinf)$ be the symmetric sequence whose $n^{th}$ space is 
  $\wedone$ bimodule maps
  \[
	\enD(\schinf)(n) \coloneqq \map_{\wedone}(\schinf(1,\smallbullet), 
	\schinf(n,\smallbullet) ).
  \]
\end{definition}

\begin{lemma}
  Operadic composition in $\schinf$ induces the structure of an operad on the 
  symmetric sequence $\enD(\schinf)$.
\end{lemma}
\proof 
  Define an operad $\mathcal E\in\Operads(\Coll{\topl})$ by setting
  \[
	\mathcal E(n,m) = \map_{\Coll{\topl}}\big(\enD(\schinf)(\smallbullet; n), 
	\,\, \schinf(\smallbullet, m)\big).
  \]
  Let $n_0 + \cdots + n_m = n$ and $\ell_1 + \cdots + \ell_m=\ell$.  The 
  monoidal structure $\otimes$ on symmetric sequences gives a map from 
  $\mathcal E(n_0,m) \times \prod_{i=1}^m \mathcal E(n_i, \ell_i)$ to
  \begin{equation}\label{equation:mapColl}
	\map_{\Coll{\topl}}\Bigg(\bigotimes_{i=0}^m \enD(\schinf)(\smallbullet; 
	n_i), \,\,\, \schinf(\smallbullet, m) \otimes \Big(\bigotimes_{i=1}^m 
	\schinf(\smallbullet,\ell_i)\Big)\Bigg).
  \end{equation}
  Now push forward from \eqref{equation:mapColl} via the operad structure on 
  $\schinf$,
  \[
	\schinf(\smallbullet,m) \otimes \Big(\bigotimes_{i=1}^m 
	\schinf(\smallbullet, \ell_i)\Big) \to \schinf(\smallbullet, \ell),
  \]
  and pull back from \eqref{equation:mapColl} by
  \begin{equation}
	\label{equation:end(schinf)(bullet;n)-to-bigotimes-end(schinf)}
		\enD(\schinf)(\smallbullet; n) \to \bigotimes_{i=0}^m \enD(\schinf) 
		(\smallbullet; n_i ).
	\end{equation}
	This defines
	\[
		\mathcal E(n_0,m) \times \prod_{i=1}^m \mathcal E(n_i, \ell_i) \to 
		\mathcal E(n,\ell).
	\]
	The morphism in 
	\eqref{equation:end(schinf)(bullet;n)-to-bigotimes-end(schinf)} comes from 
	the sequence of maps
  \begin{align*}
	\enD(\schinf)(k ; n)
	& = \coprod_{f \colon \ul k \to \ul n} \,\, \prod_{i=1}^n 
	\enD(\schinf)(f\inv(i)) \\
	& = \coprod_{f \colon \ul k \to \ul n} \,\, \prod_{j \in \ul 0 \sqcup \ul 
m} \,\, \prod_{i \in 	g\inv(j)} \enD(\schinf)(f\inv(i)) \\
	& \to \coprod_{\tilde f \colon \ul k \to \ul 0 \sqcup \ul m }\,\, \prod_{j 
	\in \ul 0 \sqcup \ul m} \enD(\schinf)(\tilde f \inv(j); n_j).
  \end{align*}
  The first equality holds  by definition, the second is a regrouping.
  The decomposition $\sum_{i=0}^m n_i = n$ defines a map $g \colon \ul n \to 
  \ul 0 \sqcup \ul m$ where $\abs{g\inv(i)}= n_i$.
  The third map sends the component corresponding to $f \colon \ul k \to \ul n$ 
  to the component corresponding to $f g \colon \ul k \to \ul 0 \sqcup \ul m$.
  For each $n,m$ there is a map,
  \begin{equation}\label{equation:schinf-to-mathcalE(n,m)}
	\schinf(n,m) \to \mathcal E (n,m).
  \end{equation}
  When $n = 0, 1$, the map \eqref{equation:schinf-to-mathcalE(n,m)} is canonical.  
  Restricting to degrees 0 and 1 gives a map of $\Operads(\Collleqone{\topl})$ 
  operads $\schone \to \mathcal E^1$, where $\mathcal E^1$ is the degree 0-1 
  part of $\mathcal E$.  Since $\schinf$ is freely generated by its degree 0 
  and 1 pieces, we get \eqref{equation:schinf-to-mathcalE(n,m)} for all $n$, 
  assembling into a map of operads in $\Operads(\Coll{\topl})$.  This 
  guarantees that \eqref{equation:schinfenD-schinf} can be used to define an 
  operadic composition law on $\enD(\schinf)$,
  \begin{equation}
	\label{equation:schinfenD-schinf}
	\schinf(n,m) \times \enD(\schinf) (k_1) \times \cdots \times 
	\enD(\schinf)(k_n) \to \schinf(k,m).
		\rlap{\hspace{0.95em}\qedsymbol} 
  \end{equation}

\begin{definition}
\label{definition:schinf-rtimes-O}
  Let $\calO$ be any 1-colored topological operad, and let $\rho \colon \calO 
  \to \enD(\schinf)$ be a map of operads.  Define the $K$-colored operad 
  $\schinf \rtimes_{\rho} \calO$ by setting
  \[
  (\schinf \rtimes_{\rho} \calO)^\h(n,m) = \schinf(n,m)
  \quad \quad
  (\schinf \rtimes_{\rho} \calO)^\e(n,m) = \begin{cases}
	\calO(n) & m = 0 \\ \emptyset & m >0
  \end{cases}
  \]
  Composition in $\schinf \rtimes_{\rho} \calO$ uses composition in $\calO$, 
  composition in $\schinf$, and the action of $\calO$ on $\schinf$ defined by 
  $\rho$,
  \begin{align*}
	\schinf(n,m) \times \calO(k;n) & \to \schinf(n,m) \times \enD(\schinf)(k;n) 
	  \to \schinf(k,m),
  \end{align*}
  where the right arrow above is the one in \eqref{equation:schinfenD-schinf}.
\end{definition}
\begin{lemma}
\label{lemma:H-is-O-algebra}
  Let $\calO$ be a topological operad and let $\rho \colon \calO \to 
  \enD(\schinf)$ be a map of operads, then the $\schinf$ structure on $(H, A)$ 
  naturally extends to a $\schinf \rtimes_\rho \calO$ structure on $(H,A)$.  In 
  particular, $H = \hoch(A)$ inherits an $\calO$ algebra structure.
\end{lemma}
\begin{proof}
  We only need to show there is a map of operads $\enD(\schinf) \to \enD(H)$ 
  compatible with the action of $\enD(\schinf)$ on $\schinf$ and the action of 
  $\schinf$ on $H$.  Indeed, the map
  \begin{equation}
  \label{equation:enDH-H}
  \enD(\schinf)(n) \otimes H^{\otimes n} \to H
  \end{equation}
  is adjoint to the maps, for all $m \geq 0$,
  \begin{align*}
  \schinf(1,m) \otimes \enD(\schinf)(n) \otimes H^{\otimes n} \otimes 
  A^{\otimes m}  \to
  \schinf(n,m) \otimes H^{\otimes n} \otimes A^{\otimes m}  \to
  A,
  \end{align*}
  where the first arrow is \eqref{equation:schinfenD-schinf} and the second arrow 
  is the $\schinf$ structure on $(H,A)$.  To check that \eqref{equation:enDH-H} 
  is compatible with composition in $\enD(\schinf)$ observe that there are two 
  morphisms of operads in $\Op{}{}(\Coll{\topl})$,
  \[
  \schinf \rightrightarrows \enD^\h(\enD(\schinf)(H), A),
  \]
  where $\enD(\schinf)(H)$ is the free $\enD(\schinf)$ algebra generated by 
  $H$.  One of the arrows uses the action of $\enD(\schinf)$ on $\schinf$, 
  while the other uses the map $\enD(\schinf)(H) \to H$ defined by 
  \eqref{equation:enDH-H}.  To check that these arrows agree, we only need to 
  check that they agree out of $\schinf(n,m)$ when $n = 0, 1$.  This is because 
  $\schinf$ is freely generated in degrees 0 and 1.  When $n = 0$, the maps are 
  obviously the same.  When $n = 1$, the maps are the same by definition of the 
  $\schone$ structure on $(H,A)$.
\end{proof}

\begin{proposition}[The universal cheese theorem]
  \label{proposition:schinf-rtimes-O}
  Let $\calO$ be a topological operad and let $\rho \colon \calO \to 
  \enD(\schinf)$ be a map of operads.  Then using the induced $\schinf 
  \rtimes_{\rho} \calO$ structure on $(\hoch(A),A)$ from 
  lemma~\ref{lemma:H-is-O-algebra} gives an isomorphism of categories
  \[
  \Alg_{(\schinf\rtimes_\rho \calO)} ^A(\calC) \cong 
  \Alg_{\calO}(\calC)_{/\hoch(A)}.
  \]
\end{proposition}
\begin{proof}
	Given any $\calC$ morphism $B \to H$ we can form the following 
	diagram.
	For brevity, we have deleted appearances $\otimes$.
	\begin{equation}
	\label{diagram:BH-cube}
	\begin{tikzpicture}[baseline=(current bounding box.center), xscale=5.5, 
			yscale=2, text height = 1.5ex, text depth = 0.25ex,
	cross line/.style={preaction={draw=white, -, line width=6pt}}]
	\node (010) at (0,1,0) {$\schinf(1,m)\calO(n)B^nA^m$};
	\node (110) at (1,1,0) {$\schinf(n,m)B^nA^m\hspace{-2em}$};
	\node (000) at (0,0,0) {$\schinf(1,m)BA^m$};
	\node (100) at (1,0,0) {$A$};
	\node (011) at (0,1,1) {$\schinf(1,m)\calO(n)H^nA^m$};
	\node (111) at (1,1,1) {$\schinf(n,m)H^nA^m$};
	\node (001) at (0,0,1) {$\schinf(1,m)HA^m$};
	\node (101) at (1,0,1) {$A$};
	\path[->]
	(010) edge (110) edge (000);
	\path[<-]
	(100) edge (110) edge (000);
	\path[->]
	(000) edge (001)
  (010) edge (011)
  (100) edge node[right, map name] {$=$} (101)
  (110) edge (111);
	\path[->]
	(011) edge[cross line] (111) edge[cross line] (001);
	\path[<-]
	(101) edge[cross line] (111) edge[cross line] (001);
 \end{tikzpicture}
\end{equation}
  Let $(B,A)$ be a $\schinf\rtimes_\rho\calO$ algebra extending the given 
  $\wedone$ structure on $A$, then by 
  corollary~\ref{corollary:AlgAschinf(C)-iso-C/Hoch(A)} we get a $\calC$ 
  morphism $B \to H= \hoch(A)$ making the right face of the cube 
  \eqref{diagram:BH-cube} commute.  The front face commutes by 
  lemma~\ref{lemma:H-is-O-algebra}.  The back face commutes by assumption. The 
  bottom face commutes by definition, and the top face commutes trivially.  
  This implies that, after composition with the maps whose codomain is $A$, the 
  left face of the cube commutes.  By adjointness, the two maps $\calO(B) 
  \rightrightarrows H$ agree, implying that $B \to H$ is indeed an $\calO$ 
  algebra morphism.

  On the other hand, given an $\calO$ algebra $B$ together with an $\calO$ 
  algebra map $B\to H$, we get an $\schinf$ structure on $(B,A)$ from the 
  underlying $\calC$ morphism.  We only need to check that the $\calO$ 
  structure on $B$ and the $\schinf$ structure on $(B,A)$ are compatible via 
  $\rho$.  Indeed, since $\schinf$ is freely generated in degrees 0 and 1, it 
  is enough to check that the back face of the cube commutes.  But this holds 
  because all other faces commute.  Most importantly, the left face commutes 
  because $B \to H$ is an $\calO$ algebra map.

  It is easy to see that each of these constructions are natural in $B$ and are 
  inverse to one another.
\end{proof}

\section{The homotopy \texorpdfstring{$E_d$}{Ed} structure on 
\texorpdfstring{$\hoch(A)$}{Hoch(A)}}
\label{section:right-Ed-action}
\addtocounter{subsection}{1}

In light of proposition~\ref{proposition:schinf-rtimes-O}, to prove the swiss 
cheese theorem, \ref{theorem:precise-statement-of-main-theorem}, we need to 
construct $\E \simeq E_d$ and
an operad morphism $\E \to \enD(\schinf)$ in such a way that the
corresponding $K$-colored operad $\schinf \rtimes \E$ is equivalent to
$\swcheese$.  While $\schinf$ has no obvious action of $E_d$, it is equivalent 
to something that does have an $E_d$ action.  The following theorem is proven 
in section \ref{section:the equivalence schinf to sch}.
\begin{theorem} \label{theorem:acyclic-cofibration}
  The natural map $\schinf\rightarrow\sch$ is an acyclic cofibration of operads 
  in $\Operads(\Coll{\topl})$.
\end{theorem}
In this section, we define the precise sense in
which theorem~\ref{theorem:acyclic-cofibration} gives us our $E_d$ action on 
$\schinf$ up to homotopy.  First, we have a lift $p$ in the following diagram,
\[
 \begin{tikzpicture}
  \matrix (m) [matrix of math nodes, row sep={3ex}, column sep={6ex}, text 
  height=1.5ex, text depth=0.25ex]
{ \schinf & \schinf \\  \sch &   \\};
\path[->]
 (m-1-1) edge node [auto, font=\scriptsize] {$\id$} (m-1-2)
 (m-1-1) edge node [left, font=\scriptsize] {$\iota$}(m-2-1);
\path[->,dashed]
 (m-2-1) edge node[right, below=2pt, font=\scriptsize] {$p$}  (m-1-2);
 \end{tikzpicture}
\]
We know that $\sch$ is cofibrant since it is obtained as the $\W$
construction applied to a $\Sigma$-cofibrant, well-pointed operad
$\swcheese^\h$ so it fits into the context covered by Berger-Moerdijk 
\cite{berger-moerdijk-w, berger-moerdijk-coloured}. In Spitzweck 
\cite{spitzweck-2001} we see that the corner axiom (or Quillen's SM7)
for monoidal  model categories tensored over topological spaces applies to 
categories of operads in topological spaces.  Thus
 we have an acyclic fibration \[
 \map(\sch, \sch) \stackrel{\iota^*}{\rightarrow} \map(\schinf, \sch)
\]
given by pre-composing with $\iota$.
 Since both $\iota p$ and $\id$ live over $\iota$, they must be homotopic.
   Let $h \colon [0,\infty]\otimes \sch  \to \sch$ be a homotopy with
   $h_0=\id$ and $h_\infty = \iota p$.

We will use this $h$ to define a homotopy right $E_d$ module structure on
$\schinf$.  For this we will use a homotopy equivalent version of $E_d$
which sits inside the $W$ construction. For simplicity we denote it by $\E$.
First we define the category $LE_d$.
This category is not monoidal, but will be used to build $\E$.  The letter
$L$ stands for \emph{level trees}. The objects of the topological category
$LE_d$ are finite sets, and the morphism space $LE_d(n,n') $ is defined to
be a quotient of
\[
\coprod_{\substack{k\geq 0\\n_1, \ldots, n_k}}
E_d(n_1,n') \times E_d(n_2, n_1) \times \cdots\times E_d(n, n_k) \times
[0,\infty]^{k}
\]
 A point of the space above is given by a sequence
$\alpha_i \in E_d(n_{i}, n_{i-1})$ for $1 \leq i \leq k+1$ and $t_i \in
[0,\infty]$ for $1 \leq i \leq k$.  For convenience of notation, we  set
$n_0 = n'$, $n_{k+1} = n$, $t_0 = \infty$, and $t_{k+1} = \infty$.  We
impose the following relations
\begin{relations}
\label{relations:LEd}
If
$t_i = 0$, then we can delete $t_i$ and replace $(\ldots, \alpha_i,
\alpha_{i+1}, \ldots)$ by the composition $(\ldots, \alpha_i \circ
\alpha_{i+1}, \ldots)$. If $n_{i} = n_{i-1}$ and
$\alpha_i$ is the identity, and $t_{i-1} = \infty = t_i$, then we can delete
$\alpha_i$ from the sequence and delete $t_i$ from the sequence.
\end{relations}
\begin{remark} In the $W$ construction, we could always delete the
appearance of an identity and sum the lengths of the surrounding edges.  We
do not allow that here since we do not have $h_{s+t} = h_s \circ h_t$.
\end{remark}
Composition in the category $LE_d$ is given by concatenating sequences,
setting the new coordinate in the factor $[0,\infty]$ between the two
sequences to be $\infty$.

We can
use the action of $E_d$  on $\sch$ as well as the maps $h_t, p, \iota$ to
define
\begin{equation}
\label{equation:LEd-map(schinf,schinf)}
 LE_d(n,n') \to \map( \schinf(n',m),\schinf(n,m) ).
\end{equation}
To do this, represent $\alpha \in LE_d(n,n')$  with a sequence $n'=n_0, n_1, 
\ldots, n_{k}, n_{k+1} =n$ together with $\alpha_i \in  E_d( n_{i}, n_{i-1})$ 
for $1 \leq i \leq k+1$ and $t_i \in [0,\infty]$ for $1 \leq i \leq k$. This 
gives a chain of maps \begin{multline}
\label{right-Ed-action-chain}
 \schinf(n',m) \xrightarrow{\iota} \sch(n_0,m) \xrightarrow{\alpha_1}
 \sch(n_1,m) \xrightarrow{h_{t_1}} \sch(n_1,m) \xrightarrow{\alpha_2} \cdots
  \\
\cdots\xrightarrow{h_{t_{k}}} \sch(n_{k},m) \xrightarrow{\alpha_{k+1}}
 \sch(n_{k+1},m) \xrightarrow{p} \schinf(n,m).
\end{multline}
The maps $\sch(n_i,m) \xrightarrow{\alpha_{i}} \sch(n_{i+1},m)$ are defined
by the action of $E_d$ on $\sch$,
\[ \sch(n_i,m) \times E_d( n_{i+1},  n_i)
\rightarrow \sch(n_{i+1},m). \]
Let us check that the relations~\ref{relations:LEd} in $LE_d$ are satisfied and 
that composition in $E_d$ corresponds to composition of maps of $\schinf$.  
Suppose $t_i = 0$ for some $i$.  Then $h_0 = \id$ so our chain of arrows 
contains
\[
 \sch(n_{i-1},m) \xrightarrow{\alpha_{i}} \sch(n_i,m) 
 \xrightarrow{\alpha_{i+1}} \sch(n_{i+1},m).
\]
The composition of these two is equal to the map given by $\alpha_i
\alpha_{i+1}\in  E_d(\underline n_{i+1}, \underline n_{i-1})$.  This is
because $\wcheeseh(-,m)$ is a right $E_d$ module.

If $n_{i} = n_{i-1}$, $\alpha_i$ is identity, and $t_{i-1}=t_i = \infty$,
then the composition $h_{t_{i-1}} \circ \alpha_i \circ h_{t_i}$ is equal to
$h_{t_{i-1}} = h_\infty$, so we are justified in deleting $\alpha_i$ and
$t_i$ from the sequence.

Now suppose we have some $t_i = \infty$, so that $\alpha \in LE_d(n,n')$
decomposes as $\beta_1\beta_2$ for some $\beta_1 \in  \E(n_i,n')$ and
$\beta_2 \in \E(n, n_{i})$.  The chain of compositions defining the action
of $\alpha$ from $\schinf(1,m)$ to $\schinf(n,m)$ contains the following
segment.
\[
 \cdots \sch(n_{i},m) \xrightarrow{h_{\infty}} \sch(n_{i},m) \rightarrow \cdots
\]
The composite of the actions of $\beta_1$ and $\beta_2$ is computed by
joining the chains for $\beta_1$ and for $\beta_2 $.  This joined chain
agrees with the chain for $\alpha$ except for the segment above, which is
replaced with the segment
\[
 \cdots \rightarrow \sch(n_{i},m) \xrightarrow{p} \schinf(n_{i},m) 
 \xrightarrow{\iota} \sch(n_{i},m) \rightarrow \cdots
\]
Since $h_\infty = \iota p$, these chains of maps have the same composition.

The maps \eqref{equation:LEd-map(schinf,schinf)} define a functor
\begin{equation}
\label{equation:LEd-enD(schinf)}
LE_d \to \enD(\schinf).
\end{equation}
 There is no obvious operad structure on $LE_d$ so we take
the smallest operad containing $LE_d$. More precisely, 
\eqref{equation:LEd-enD(schinf)} is a morphism of topological categories whose 
objects are finite sets.  There is a forgetful functor from operads to the 
category of such topological categories.  The operad $\E$ is defined to be the 
result of applying the left adjoint of this forgetful functor to the category 
$LE_d$.
\begin{definition}
  Let $F(LE_d)$ be the free one-colored operad generated by the symmetric 
  sequence $n
  \mapsto LE_d(n,1)$.  For each $n, n'\geq 0$ let $\E(n,n')$ be the topological
  space given by the coequalizer
  \[
	\coprod_{n''} LE_d(n'',n') \times LE_d(n,n'') \rightrightarrows 
	F(LE_d)(n,n') \xrightarrow{\scriptstyle\text{eq}} \E(n,n'),
  \]
  where the two maps are given by composition in either $F(LE_d)$ or $LE_d$ and 
  the inclusion of $LE_d$ into $F(LE_d)$.
\end{definition}

\begin{lemma}
\label{lemma:E-simeq-Ed}
  The category $\E$ is an operad and is equivalent to $E_d$.
\end{lemma}
\begin{proof}
	Given a tree with its internal edges labeled by lengths $[0,\infty]$, call 
	it a \emph{level tree} if edges equidistant from the root vertex have the 
	same length.  Every morphism in $LE_d(n,1)$ can be represented by a level 
	tree with vertices labeled by $E_d$.  We can represent a point of $F(LE_d)$ 
	with a tree whose vertices are labeled by level trees in $LE_d$.  The 
	relation defining $F(LE_d) \to \E$ allows us to break up a level tree with 
	at least one level of length $\infty$ into several level trees all of whose 
	levels have finite length. We conclude that $\E$ consists of trees labeled 
	by $E_d$ on the vertices, and $[0,\infty]$ on the internal edges, 
	satisfying the condition that every maximal finite subtree is level.

There is an operad morphism $\E \to E_d$ which collapses all edge lengths to 0.  
On the level of symmetric sequences, there is a homotopy inverse $E_d \to \E$.  
The homotopy $g_t \colon \E \to \E$ first collapses lengths of the edges 
furthest from the root to zero.	This preserves the condition that every maximal 
finite subtree is level.  Continuing in this way, we collapse all edge lengths 
to zero.
\end{proof}

The adjoint to \eqref{equation:LEd-enD(schinf)} is an operad morphism $\E \to 
\enD(\schinf)$, which by definition~\ref{definition:schinf-rtimes-O} we can use 
to define the $K$-colored operad $\schinf \rtimes \E$.

\begin{lemma}\label{lemma:schinf-rtimes-E-simeq-swcheese}
The $\{ \e,\h \}$-colored operad $\schinf \rtimes \E$ is weakly equivalent to 
the swiss cheese operad.
\end{lemma}
\begin{proof}
  First, note that $\swcheese$ is equivalent to the semi-direct
product of $\sch$ and $  \wed$ where the action of $ \wed$ factors through
the map $   \wed \to    E_d$ which sends all lengths of internal edges to
zero.  This is because the map $\sch \to \swcheese^\h$ which collapses trees
is a weak equivalence and respects the action of $\wed$.

The action of $\E$ on $\schinf$ can be extended to an action on
all of $\sch$.  The sequence $\alpha_1, t_1, \cdots, t_k , \alpha_{k+1}$ acts
via the composition \begin{multline} \label{wed-action-on-sch} \sch(n_0, m)
\xrightarrow{\alpha_1} \sch(n_1,m) \xrightarrow{h_{t_1}} \sch(n_1,m) \to
\cdots \to \sch(n_{l-1}, m)\\ \xrightarrow{h_{t_l}} \sch(n_{l-1}, m)
\xrightarrow{\alpha_l} \sch(n_l, m) \xrightarrow{h_\infty} \sch(n_l,m).
\end{multline}
Define for each $s \in [0,\infty]$ a homotopy $h^{[0,s]}\colon 
[0,\infty]\otimes \sch
 \to \sch$ by setting $h^{[0,s]}_t = h_{\min(s,t)}$.  We have $h^{[0,s]}_0 = 
 \id$ and $\h^{[0,s]}_{\infty} = h_s$, therefore we can define an action of 
 $\E$ on $[0,\infty]_s \otimes\sch$ by replacing $h_{t_i}$ in 
 \eqref{wed-action-on-sch} with
$h^{[0,s]}_{t_i}$. Then, when $s=0$ each $h^{[0,0]}_{t_i}$ is the identity,
so the action factors through the map $\E \to E_d$ collapsing all edges to 0.  
When $s=\infty$ we have $h^{[0,\infty]}_{t_i} =
h_{t_i}$ so the action of $\E$ on $\sch$ is \ref{wed-action-on-sch}.  Thus we 
have a diagram of
equivalences \[ \swcheese \leftarrow \sch\rtimes_{s=0} \E
\rightarrow ( [0,1] \otimes\sch)\rtimes  \E \leftarrow \sch
\rtimes_{s=\infty}      \E \rightarrow \schinf \rtimes \E, \]
where the map on the left collapses all edge lengths to 0.
\end{proof}

\section{The equivalence \texorpdfstring{$\schinf \to \sch$}{SChoo to SCh}}
\label{section:the equivalence schinf to sch}
This section is dedicated to proving theorem \ref{theorem:acyclic-cofibration}.  
The proof uses a recasting (\ref{definition:an-operad-is-a-functor}) of 
definition \ref{definition:colored-operad} which is equivalent when considering 
operads in $\topl$ \cite{getzler-operads-revisited}.  First we set the context for this new 
definition, then we prove that $\schinf \to \sch$ is a cofibration. Finally, we 
show that it is a weak equivalence.

\subsection{The category of Forests}
\label{section:forests}

The following definition is an amalgamation of those found in 
\cite{kontsevich-soibelman-deformations}, \cite{costello-a-infinity} and 
\cite{getzler-operads-revisited}.

\begin{definition}\label{definition:category-of-forests}
Fix a set $K$.  A \emph{$K$-colored young forest} is an uncolored map of finite 
$K$-colored  sets $x \colon I_x \to J_x$.  A \emph{$K$-colored forest} $f 
\colon x \to y$ is a color-preserving isomorphism $f\colon I_y\sqcup J_x \to 
J_y \sqcup I_x$ such that for every $i \in I_x$ there is a $k\geq 0$ such that 
$(f\circ x)^k(i) \in J_y$ and for every $i \in I_y$ there is an $\ell \geq 0$ 
such that $(f\circ x)^\ell \circ f(i) = y(i)$.
\end{definition}

\begin{definition}
	\label{definition:edges vertices of forest}
	Given a $K$-colored forest $f \colon x \to y$ we call $V(f)\coloneqq J_x$ 
	the set of \emph{internal vertices} of $f$.  We call $\inverts(f)\coloneqq 
	I_y$ the set of \emph{input vertices} of $f$ and $\rt(f)\coloneqq J_y$ the 
	set of \emph{root vertices} of $f$. In addition, $\exEdges(f) = J_x\sqcup 
	I_y \cong J_y \sqcup I_x$ is called the set of \emph{extended edges} of $f$ 
	and $E(f) \coloneqq J_x \times_f I_x$ is the set of \emph{internal edges} 
	of $f$.
\end{definition}

\begin{definition}
	\label{definition:composition of forests}
	If $g \colon x \to y$ and $f \colon y \to z$ are forests, we can define a  
	composite forest $fg \colon x \to z$.  We use concatenation to denote this 
	composition and $\circ$ to denote composition of maps of finite sets. The 
	forest $fg$, as a map $I_z \sqcup J_x \to J_z \sqcup I_x$, is defined by 
	the following rule.  If $i \in I_z$, then there is a $k \geq 0$ and an 
	$\varepsilon \in \{0,1\}$ such that $g^\varepsilon (f\circ g)^k(i) \in J_z 
	\sqcup I_x$.  Similarly, if $i \in J_x$ there is a $k$ and $\varepsilon$ 
	such that $f^\varepsilon (g\circ f)^k (i) \in J_z \sqcup I_x$.  In 
	\cite{thomas-forests-arxiv}, it is shown that forest composition is 
	associative.  This, together with disjoint union, makes young forests the 
	objects and forests the morphisms of a symmetric monoidal category denoted 
	$\for$.
	
\end{definition}

\begin{definition}
	\label{definition:an-operad-is-a-functor}
	A $K$-colored operad $\calO$ is a strong symmetric monoidal functor 
	$(\rfor_K, \sqcup) \to (\topl,\times)$.  The category of operads is the 
	category $\fun^{\otimes}(\for_K, \topl)$ of symmetric monoidal functors and 
	natural transformations.  We denote this category by $\Op_K$ just as in 
	definition~\ref{definition:opK}.
\end{definition}

\begin{remark}
	The category of $K$-colored operads $\Op_{K}$ as defined in 
	\ref{definition:opK} is naturally isomorphic to the functor category 
	$\fun^\otimes(\for_K, \topl)$.  Indeed, given $\calO$ from 
	\ref{definition:colored-operad} as in definition~\ref{definition:opK} we 
	define $\calO(x)$ for a young forest $x$ to be $\bigotimes_{j\in J_x} 
	\calO(x\inv(j);j)$.
\end{remark}

\begin{remark}
  \label{remark:drop-K-from-notation}
  When the set of colors $K$ is understood,
  we often drop it from the notation.  Forests and young forests are always 
  $K$-colored, for some set $K$.  The category $\for_K$ will be abbreviated 
  $\for$, and $\Op_{K}$ will be denoted $\Op$.
\end{remark}

\begin{definition}
	\label{definition: fibrations in categories of symmetric monoidal functors}
	Let $(C, \otimes)$ be a symmetric monoidal category. Call an object $c \in 
	C$ \emph{indecomposable} if it cannot be written as a tensor product $c 
	\cong c_1 \otimes c_2$ for any $c_1, c_2 \in C$.  Let $\fun^{\otimes}(C, 
	\topl)$ denote the category of strong symmetric monoidal functors 
	$(C,\otimes) \to (\topl, \times)$.  We call a morphism $\psi \colon \calO 
	\to \calP$ in $\fun^\otimes(C, \topl)$ a \emph{fibration} (respectively 
	\emph{weak equivalence}) if $\psi(c) \colon \calO(c) \to \calP(c)$ is a 
	fibration (respectively weak equivalence) for every indecomposable $c \in 
	C$.  Define the class of \emph{cofibrations} in the usual manner (see Hovey  
	\cite{hovey-book}).
\end{definition}

For every symmetric monoidal category $C$ we consider in this paper we will use 
\ref{definition:
fibrations in categories of symmetric monoidal functors} to define 
cofibrations, fibrations, and weak equivalences regardless of whether or not 
these form a model structure.

\begin{remark}
	\label{remark: when a forest is a tree}
	A young forest $x$ is indecomposable (definition~\ref{definition: 
	fibrations in categories of symmetric monoidal functors}) if $J_x \simeq 
	*$. In this case we say $x$ is a young \emph{tree}.   We say a forest $f 
	\colon y \to x$ is a \emph{tree} if $x$ is a young tree.
\end{remark}

\begin{definition}
	\label{definition:model structure on operads as functors}
	For a category $C$, let $C^\times$ denote the category with the same 
	objects as those of $C$, but only invertible morphisms.  Let $C_{\indec}$ 
	denote the full subcategory (not monoidal) of $C$ consisting indecomposable 
	objects.  Let $C^\times_{\indec} = (C_{\indec})^\times$.  The functor 
	category $\topl^{C^\times_{\indec}}$ is called the category of \emph{
	$C$-symmetric sequences}, denoted $C\dash\SymSeq$.
\end{definition}

\begin{remark}
	A forest $f \colon x \to y$ is invertible if and only if $f(J_x) = J_y$ and 
	$f(I_y) = I_x$.  Thus an invertible forest gives a pair of isomorphisms of 
	$K$-colored finite sets $J_f \colon J_x \to J_y$ and $I_f \colon I_y \to 
	I_x$ which are compatible with the maps $x \colon I_x \to J_x$ and $y 
	\colon I_y \to J_y$.  We conclude that the category $\for^\times_\indec$ is 
	isomorphic to the \emph{opposite} of the category of $K$-colored finite 
	sets and bijections.   In the case $K\simeq \{*\}$ we get 
	$\for^\times_\indec \cong \Sigma$ where $\Sigma$ is as in
	definition~\ref{definition:symmetric sequences}.  Moreover, in this case 
	the category of $\for$-symmetric sequences is isomorphic to the category 
	$\topl^{\Sigma}$ of symmetric sequences.
\end{remark}

\begin{notation}
	\label{notation:FCD}
	If $C \to \for$ is any symmetric monoidal functor, we will denote the 
	category $\fun^\otimes(C,\topl)$ of strong symmetric monoidal functors  by 
	$C\dash\Op$ unless we say otherwise (for example, we do not use this 
	notation in definition~\ref{definition: Op_omega and Op_k}). We call the 
	objects of $C\dash\Op$ \emph{$C$-operads}. In all cases we consider, the 
	functor $C \to \for$ will be apparent from the category $C$ so we leave the 
	functor out of the notation.  If $C \to D \to \for$ is a pair of symmetric 
	monoidal functors we denote the forgetful functor $D\dash\Op \to C\dash\Op$ 
	by $U^D_C$, and denote its left adjoint by $F^C_D$.
\end{notation}

%
%
%

\subsection{The \texorpdfstring{$W$}{W} construction}
\label{section:w-construction}
We show how the $\W$ construction of Boardman and Vogt 
\cite{boardman-vogt-73} can be realized as a coend construction using 
$\for$.  In nice situations the $W$ construction gives a cofibrant replacement 
for an operad, as shown in \cite{berger-moerdijk-w}.  We use $[0,\infty]$ as 
our edge labels as in \cite{kontsevich-operads-motives}.

Given any young forest $z$ there is a contravariant functor $W \colon 
\rfor_{/z}\to \topl$ from the over category of $z$ to topological spaces.  For 
any object $g \colon y \to z$ of this over category, set $W(g) = \map(E(g), 
[0,\infty])$.  If $f \colon x \to y$ is a forest, the map $W(g) \to W(gf)$ is 
denoted $W_\Sigma(f)$. This map uses the sum operation on $[0,\infty]$.  This 
is an extension of $+$ on $[0,\infty)$ such that $t + \infty = \infty = \infty 
	+ t$ for all values of $t$.  Concretely, $g$ and $gf$ define maps of sets
\[
	E(g) \xrightarrow{\tilde{g}} J_z \sqcup I_x \xleftarrow{\tilde{gf}} 
	E(gf)
\]
We can turn a function $t \in W(g)$ to a function $W_\Sigma(t) \in W(gf)$ by 
pushing forward along $\tilde{g}$, then pulling back along $\tilde{gf}$.  
Pushing forward means summing over fibers, which is well-defined since all the 
sets we are considering are finite.
	

If $h \colon z \to w$ is a forest there is a morphism $W_\infty(h) \colon W(g) 
\to W(hg)$ which uses the maps $E(g) \hookrightarrow I_y \hookleftarrow E(hg)$.  
In this case we do not push forward and pull back functions.  Rather we extend 
a function $t \colon E(g) \to [0,\infty]$ to a function on $E(hg)$ by setting 
$t(\epsilon) = \infty$ if $\epsilon \not\in E(g)$.  This defines a natural 
transformation $W_\infty(h) \colon W \to W h_*$ where $h_* \colon 
\rfor^{op}_{/z} \to \rfor^{op}_{/w}$ is induced by $h \colon z \to w$.  In 
other words, diagram~\eqref{diagram:Wg-Wgf-Whgf} commutes for all composeable 
triples $h, g$, and $f$.
\begin{equation}\label{diagram:Wg-Wgf-Whgf}
	   \begin{tikzpicture}[numbered picture]
		  \matrix (m) [matrix of math nodes, row sep={2.5ex}, column 
			  sep={2.5em},
			 text height=1.5ex, text depth=0.25ex]
		  {
		   W(g) & W(gf)
		  \\
		   W(hg) & W(hgf)
		  \\};
		  \path[->, every node/.style={map name}]
		   (m-1-1) edge node[above] {$W_\Sigma (f)$} (m-1-2) edge node[left] 
		   {$W_\infty(h)$} (m-2-1)
		   (m-1-2) edge node[right] {$W_\infty (h)$} (m-2-2)
		   (m-2-1) edge node[above] {$W_\Sigma (f)$} (m-2-2);
	   \end{tikzpicture}
	\end{equation}

Consider an operad $\calO$ as a collection of functors $\calO_z \colon 
\rfor_{/z} \to \topl$ by setting $\calO_z(g) \coloneqq \calO(y)$ for $g\colon y 
\to z$.  For a young forest $z$ the topological space $W\calO(z)$ is the coend
\begin{equation}
	\label{equation:WO(z)-as-coend}
	\W \calO (z) = W\otimes_{\rfor_{/z}}\calO_z = \left(\coprod_{g \colon y \to 
	z} W(g) \times \calO(y) \right)\modsim
\end{equation}
where $(W_{\Sigma}(f)t, \alpha)\sim(t, \calO(f)\alpha)$ for every $f \colon x 
\to y$, $t \in W(g)$ and $\alpha \in \calO(x)$.
Now a forest $h \colon z \to w$ with the natural transformations above gives us 
a map
\begin{equation}
  \label{equation:WO(z)-to-WO(w)}
  W\calO(z) = W \!\otimes_{\rfor_{/z}} \!\calO_z \xrightarrow{W_\infty(h) 
  \otimes \id} W h_* \!\otimes_{\rfor_{/z}}\!\calO_z \to 
  W\!\otimes_{\rfor_{/w}} \!\calO_w = W\calO(w).
\end{equation}
This defines $W\calO$ as a functor $\rfor \to \topl$.  This functor is 
symmetric monoidal, so $\W\calO$ is a $K$-colored operad.

\begin{remark}
	\label{remark: labeled forests}
	The space $W\calO(z)$ consists of labeled forests $f\colon y \to z$.  A 
	labeled forest is one whose internal edges each have a length in 
	$[0,\infty]$ and whose vertices $j \in J_x$ have a corresponding label in 
	$\calO(x\inv(j); j)$.
\end{remark}

In the sequel, we will define several variants on the category $\for$.  Each of 
these variants admits a functor to $\for$ and we want a corresponding $W$ 
construction for each.
\begin{definition}
	\label{definition:WvarO}
	Suppose $\calO$ is an operad $\for \to \topl$ and $C$ is any symmetric 
	monoidal category equipped with a symmetric monoidal functor $G\colon C \to 
	\for$. From this data we can construct an operad $W_{C}\calO \in C\dash\Op$ 
	using formulas analogous to
	\eqref{equation:WO(z)-as-coend} and \eqref{equation:WO(z)-to-WO(w)}.  
	Specifically,  let $W_{C}$ be the composite $W\circ G_{c} \colon 
	(C_{/c})^{op} \to (\for_{/G(c)})^{op} \to \topl$ where $c \in C$.  We 
	define
	\begin{equation}
		\label{equation:WCO(c)-as-coend}
		W_{C} \calO (c) = W_C\otimes_{C_{/c}}(G^*\calO)_c = \left(\coprod_{f 
		\colon b \to c} W(G(f)) \times \calO(G(y)) \right)\modsim
	\end{equation}
	We define $W_C\calO(c) \to W_C\calO(d)$ for a $C$-morphism $g \colon c \to 
	d$ just as in \eqref{equation:WO(z)-to-WO(w)}.
\end{definition}

\begin{example}
	\label{example:sch as w construction}
	Define the full subcategory $D \hookrightarrow \for$ to be given by those  
	young forests $x$ where $J_x$ has only color $\h$.  Restricting the swiss 
	cheese operad to this full subcategory gives the operad $\sch$ as 
	$W_{D}\swcheese$.  To get $\schone$, we use $C \hookrightarrow \for$, the 
	full subcategory of young forests $x$ where $J_x$ has color $\h$ and for 
	each $j \in J_x$, there is at most one element of $x\inv(j)$ of color $\e$.

	It should be clear that $D \dash\Op \simeq \Op(\Coll{\topl})$.  Indeed, the 
	forgetful functor $\Op_{K} \to \Op(\Coll{\topl})$ from 
	definition~\ref{definition:opK} is given by pulling back along the 
	symmetric monoidal functor $D \to \for_K$.    Finally, note that $\schinf$ 
	as defined in definition~\ref{definition:schinf} is
	$F^{C}_{D} \schone$, where $F^{C}_{D}$ is defined in notation 
	\ref{notation:FCD}
\end{example}

\subsection{Weighted Forests}
\label{subsection:weighted forests}

We need to define the category of weighted forests to prove the following half 
of theorem \ref{theorem:acyclic-cofibration}.
\begin{theorem}
	\label{theorem:schinf to sch cofibration}
	The natural map $\schinf \to \sch$ is a cofibration in $\Op(\Coll{\topl})$.
\end{theorem}
The proof follows closely the work of Berger-Moerdijk \cite{berger-moerdijk-w}.  
From an operad $\calO$ they construct an increasing chain of symmetric 
sequences
\begin{equation}
	\label{equation:W0O to W1O etc}
	W_0\calO \to W_0^+ \calO \to W_1 \calO \to W_1^+ \calO \to W_2 \calO \to 
	\cdots
\end{equation}
The symmetric sequence $W_0 \calO$ is just the underlying symmetric sequence of 
$\calO$.  If $\calO$ is cofibrant as a symmetric sequence $\calO$ is said to be 
\emph{$\Sigma$-cofibrant}.  If the operadic unit maps of $\calO$ are 
cofibrations, $\calO$ is said to be \emph{well-pointed}.  In the case $\calO$ 
is well-pointed and $\Sigma$-cofibrant, Berger-Moerdijk show that $W\calO 
\coloneqq \colim_k W_k\calO$ is a cofibrant replacement of $\calO$ as an 
operad.  In the course of the proof, they show that $W_k \calO$ is a 
\emph{$k$-operad}, which is a partial operad in a certain sense.  This partial 
operad structure will be encoded
here in the category $\Op_k$ from definition~\ref{definition: Op_omega and 
Op_k}.  Each $W_k^+ \calO$ is an operad in $\Op_{k+1}$, and in our context, is 
given by the left adjoint to a forgetful functor $\Op_{k+1} \to \Op_{k}$ 
applied to $W_k\calO$.  Concretely, the points of $W_k\calO$ are given by trees 
with at most $k$ internal edges whose vertices are labeled by $\calO$ and whose 
internal edges are labeled by $[0,\infty]$.  This section mimics this work of 
Berger-Moerdijk to prove theorem~\ref{theorem:schinf to sch cofibration}.

	

\begin{definition}
	\label{definition:(f|x)}
	Let $f \colon x \to y$ be a forest.  Let $(f|x)$ denote the endomorphism of 
	$I_y \sqcup J_y \sqcup I_x \sqcup J_x$ which is $f$ on $I_y \sqcup J_x$, $x$ 
	on $I_x$, and the identity on $J_y$.  By assumption that $f$ is a forest, 
	there is a $k\geq 0$ such that for every $i$ we have $(f|x)^k(i) \in J_y$.  
	Let $[f|x]$ denote $(f|x)^\infty$.
\end{definition}
\begin{definition}
	\label{definition:weighted-forests}
	For $I \subset I_y \sqcup I_x \sqcup J_x$ and $j \in J_y$ let $I(j)$ denote 
	$I\cap [f|x]\inv(j)$, the set of elements of $I$ living over $j$.   A 
	\emph{weighted young forest} is a pair $(x, \omega_x)$ where $x$ is a young 
	forest and $\omega_x \colon J_x \to \mathbb Z_{\geq 0}$ is any function, 
	called the \emph{weight of $x$}.    A \emph{weighted forest} $f \colon (x, 
	\omega_x) \to (y, \omega_y)$ is a forest $f \colon x \to y$ such that for 
	all $j \in J_y$,
	\begin{equation}
		\label{equation:omega-y-in-terms-of-f-and-omega-x}
		\omega_y(j) \geq \#E(f) + \sum_{{i \in J_x(j)}} \!\!\!\omega_x(i).
	\end{equation}
	If $g \colon (y,\omega_y) \to (z, \omega_z)$ is a weighted forest.
	then one can show that
	$gf \colon x \to z$ defines a weighted forest $gf \colon (x, \omega_x) \to (z, 
	\omega_z)$.

\end{definition}

\begin{definition}
	\label{definition: Op_omega and Op_k}
  Disjoint union of forests extends to disjoint union of weighted forests. Let 
  $\rfor_\omega$ denote the symmetric monoidal category of weighted forests.
  For each $k \geq 0$, let $\rfor_{k}$ denote the full subcategory of 
  $\rfor_\omega$ generated by objects of the form $(x,\omega_x)$ such that 
  $\omega_x(j) \leq k$ for every $j \in J_x$.  Let $\Op{k}$ denote the category 
  of \emph{($K$-colored) weight $k$ operads}, which are strong symmetric 
  monoidal functors $\rfor_k \to \topl$.  Note that there is a symmetric 
  monoidal functor $\for_\omega \to \for$ which forgets the weights.

  
\end{definition}

\begin{remark}
	\label{remark: Sigma cofibrance and for0}
	If a forest $f \colon x \to y$ has no internal edges then $f$ is a disjoint 
	union of isomorphisms and maps of the form $[\emptyset \to \emptyset] \to 
	[\{\kappa\} \to \{\kappa\}]$, where $\kappa \in K$.  The restriction of an 
	operad $\calO$  to $\for_0$ remembers only
	\begin{itemize}
		\item the spaces $\calO(I;\kappa)$ for each $K$-colored set $I$ and 
			each color $\kappa$, together with the right $\aut(I)$ action on 
			$\calO(I;\kappa)$, and
		\item the operadic unit maps $* \to \calO(\kappa;\kappa)$ for each 
			color $\kappa \in K$.
	\end{itemize}
	Thus $\calO$ is cofibrant  (definition~\ref{definition: 
	fibrations in categories of symmetric monoidal functors}) as an object of 
	$\fun^{\otimes}(\for_0, \topl)$ if and only if it is well-pointed and 
	$\Sigma$-cofibrant as in \cite[section 3]{berger-moerdijk-w}.
\end{remark}

\begin{proof}[Proof of theorem \ref{theorem:schinf to sch cofibration}]
	In example \ref{example:sch as w construction} we constructed categories 
	$C$ and $D$ such that $\schone = W_C\swcheese$, $\sch = W_D\swcheese$ and 
	the map $\schinf \to \sch$ is $F_D^C W_C\swcheese \to W_D \swcheese$.  
	Lemma~\ref{lemma: FWCO to WDO cofibration} applies since $C$ is a full 
	subcategory of $D$ which is a full subcategory of $\for$,
	and $\swcheese \colon \for \to \topl$ is cofibrant as a functor $\for_0 \to 
	\topl$ (see remark~\ref{remark: Sigma cofibrance and for0}).
\end{proof}




\begin{definition}
	\label{definition: C omega}
	Let $C$ be any full subcategory of $\for$.  Define $C_\omega$ to be the 
	full subcategory of $\for_\omega$ given by young weighted forests $(x, 
	\omega_x)$ such that $x \in C$.  The functor $\for_\omega \to \for$ induces 
	a functor $C_\omega \to C$.
\end{definition}
If $C \hookrightarrow \for$ is a full subcategory, there is a map of $C_\omega$ 
operads $W_{C_\omega} \to U^{C}_{C_\omega} W_{C}\calO$.  Concretely, for each 
weighted young forest $(x,\omega_x)$ we have a map
\begin{equation}
	\label{equation: WComegaO to UWCO before relation}
	\coprod_{f \colon (y, \omega_y) \to (x, \omega_x)} W(f) \times \calO(y) \to
	\coprod_{f \colon y \to x} W(f) \times \calO(y),
\end{equation}
defined in the obvious way.  This descends to give a morphism
\begin{equation}
	\label{equation: WComegaO to UWCO}
	W_{C_\omega}\calO \to U^{C}_{C_\omega} W_C \calO
\end{equation}
\begin{lemma}
	\label{lemma:Fomega Womega is W}
	For any full subcategory $C \hookrightarrow \for$, the left adjoint of 
	\eqref{equation: WComegaO to UWCO},
	\[
		F_C^{C_\omega} W_{C_\omega} \calO  \to W_C \calO,
	\]
	is an isomorphism.
\end{lemma}
\begin{proof}
	The left adjoint $F^{C_\omega}_C W_{C_\omega} \calO$ can be computed at a 
	young tree $x$ as $\colim_{k \to \infty}(W_{C_\omega} \calO)(x, k)$.  For 
	each $k$, a point of $(W_{C_\omega} \calO)(x, k)$ is given by a labeled 
	tree (\ref{remark: labeled forests}) $f \colon y \to x$ with at most $k$ 
	internal edges.  Taking the colimit as $k$ goes to $\infty$, we get all 
	labeled trees over $x$, which is $W_C\calO(x)$.
\end{proof}

\begin{lemma}
	\label{lemma:when a cofibration of symmetric sequences is cofibration of 
	operads}
	Let $C$ be any symmetric monoidal full subcategory of $\for_\omega$.
	Suppose $\psi \colon \calO \to \calP$ is a morphism of operads in 
	$C\dash\Op$ which is a cofibration of $C$-symmetric sequences. Further 
	suppose that for every young forest $x \in C$ such that there is a 
	\emph{tree} $f \colon x \to y$ in $C$ with at least one internal edge we 
	have that $\psi(x) \colon \calO(x) \to \calP(x)$ is an isomorphism.  These 
	conditions imply that $\psi$ is a cofibration of $C$-operads.
\end{lemma}
\begin{proof}
	Suppose $\calQ \to \calQ'$ is an acyclic fibration of $C$-operads.  By 
	definition~\ref{definition: fibrations in categories of symmetric monoidal 
	functors} this means that for every young tree $x \in C$, $\calQ(x) \to 
	\calQ'(x)$ is an acyclic fibration of topological spaces.  Suppose we have 
	a commutative diagram of $C$-operads
	\begin{equation}
		\label{diagram:commutative diagram of C operads}
		\begin{tikzpicture}[
				numbered picture,
				node distance=6ex and 3em,
				math nodes,
				text height=2ex,
				text depth=0.75ex,
				every node/.style={on grid,inner sep=1pt},
				mn/.style={inner sep=1pt, font=\scriptsize}
			]
			\node (ul) {\calO};
			\node[right=of ul] (ur) {\calQ};
			\node[below=of ul] (ll) {\calP};
			\node[right=of ll] (lr) {\calQ'.};
			\path[->]
			(ul) edge (ur)
			(ul) edge (ll)
			(ur) edge (lr)
			(ll) edge (lr);
		\end{tikzpicture}
	\end{equation}
	By assumption, there is a lift $\rho \colon \calP \to \calQ$ \emph{of 
	$C$-symmetric sequences}.  We claim that $\rho$ is automatically a morphism 
	of $C$-operads.  To prove this claim, it is enough to show that the square 
	on the right in diagram~\eqref{diagram: P to Q an operad morphism} commutes 
	for every tree $f \colon x \to y$ in $C$.
	\begin{equation}
		\label{diagram: P to Q an operad morphism}
		\begin{tikzpicture}[
				numbered picture,
				node distance=7ex and 5em,
				math nodes,
				text height=2ex,
				text depth=0.75ex,
				every node/.style={on grid,inner sep=1pt},
				mn/.style={inner sep=1pt, font=\scriptsize}
			]
			\node (ul) {\calP(x)};
			\node[right=of ul] (ur) {\calQ(x)};
			\node[below=of ul] (ll) {\calP(y)};
			\node[right=of ll] (lr) {\calQ(y)};
			\node[left=of ul]  (ull) {\calO(x)};
			\node[below=of ull] (lll) {\calO(y)};
			\path[->]
			(ul) edge node[above, mn] {\rho(x)} (ur)
			(ul) edge node[left, mn] {\calP(f)} (ll)
			(ur) edge node[right, mn] {\calQ(f)} (lr)
			(ll) edge node[above, mn] {\rho(y)} (lr)
			(ull) edge node[above, mn] {\psi(x)} (ul)
			(ull) edge node[left, mn] {\calO(f)} (lll)
			(lll) edge node[above, mn] {\psi(y)} (ll)
			;
		\end{tikzpicture}
	\end{equation}
	Suppose $f$ is an isomorphism. Then we know the square commutes because 
	$\rho$ is a map of $C$-symmetric sequences.  Suppose $f$ is not an 
	isomorphism and has no internal edges. Since $f$ is a tree it must be of 
	the form $[\emptyset \to \emptyset] \to [\{c\} \to \{c\}]$ for some $c \in 
	K$.  In this case we have $\calO(x) = \calP(x) = \calQ(x) = *$, in 
	particular $\psi(x)$ is an isomorphism.  Thus our assumption shows that if 
	$f$ is not an isomorphism then $\psi(x)$ is an isomorphism.  We can deduce 
	that the square on the right commutes in this case from the fact that the 
	square on the left and the outer square commute.
\end{proof}

Given a full subcategory $D \hookrightarrow \for$ and a full subcategory $C 
\hookrightarrow D$, we can interpolate between $C_\omega$ and $D_\omega$ with a 
sequence of subcategories of $\for_\omega$.
For each $k \geq 0$, let $D_k$ be the full subcategory of $\for_\omega$ given 
by disjoint unions of young weighted trees $(x,\omega_x)$ where either $x \in 
C$ or $x \in D$ and $\omega_x \leq k$.  Note that $D_{-1} = C_\omega$.  We have 
left the inclusion $C\hookrightarrow D$ implicit in the notation $D_k$.  We 
have a commutative diagram of symmetric monoidal functors
\begin{equation}
	\label{diagram:sequence of full subcategories}
	\begin{tikzpicture}[
			numbered picture,
			node distance=5ex and 3em,
			math nodes,
			text height=2ex,
			text depth=0.75ex,
			every node/.style={on grid,inner sep=1pt},
			mn/.style={inner sep=1pt, font=\scriptsize}
		]
		\node (C) {C};
		\node[below=of C] (D-1) {C_\omega=D_{-1}};
		\node[right=4em of D-1] (D0) {D_0};
		\node[right=of D0] (cdots) {\cdots};
		\node[right=of cdots] (Domega) {D_\omega.};
		\node[above=of Domega] (D) {D};
		\path[->]
		(C) edge (D)
		(D-1) edge (C) edge (D0)
		(D0) edge (cdots)
		(cdots) edge (Domega)
		(Domega) edge (D);
	\end{tikzpicture}
\end{equation}
The bottom row consists of successively larger full subcategories of 
$\for_\omega$ and the top row consists of full subcategories of $\for$. 

From the symmetric monoidal functors in  \eqref{diagram:sequence of full 
subcategories} we get the categories and forgetful functors in 
\eqref{diagram:sequence of forgetful functors}.
\begin{equation}
	\label{diagram:sequence of forgetful functors}
	\begin{tikzpicture}[
			numbered picture,
			node distance=5ex and 5em,
			math nodes,
			text height=2ex,
			text depth=0.75ex,
			every node/.style={on grid,inner sep=1pt},
			mn/.style={inner sep=1pt, font=\scriptsize}
		]
		\node (C) {C\dash\Op};
		\node[below=of C] (D-1) {C_\omega\dash\Op = D_{-1}\dash\Op};
		\node[right=7em of D-1] (D0) {D_0\dash\Op};
		\node[right=of D0] (cdots) {\cdots};
		\node[right=of cdots] (Domega) {{D}_\omega\dash\Op};
		\node[above=of Domega] (D) {D\dash\Op};
		\path[<-]
		(C) edge (D)
		(D-1) edge (C) edge (D0)
		(D0) edge (cdots)
		(cdots) edge (Domega)
		(Domega) edge (D);
	\end{tikzpicture}
\end{equation}

\begin{remark}
	\label{remark: F preserves cofibrations}
	The forgetful functor $D_{k}\dash\Op \to D_{\ell}\dash \Op$ for $k \geq 
	\ell$, where $k, \ell \in \{-1, 0, 1, \ldots\} \sqcup \{\omega\}$ clearly 
	preserves fibrations (definition~\ref{definition: fibrations in categories 
	of symmetric monoidal functors}), thus the left adjoint 
	$F_{D_{k}}^{D_\ell}$ preserves cofibrations.
\end{remark}


\begin{lemma}
	\label{lemma: colim FWkO is WomegaO}
	For each $k \geq 0$ there is a natural map in $D_\omega\dash\Op$,
	\[
		F^{D_k}_{D_\omega} W_{D_k} \calO \to W_{D_\omega} \calO.
	\]
	The colimit of these maps as $k\to \infty$ is an isomorphism.
\end{lemma}
\begin{proof}
	For each young tree $x$ and each weight $\omega_x = \ell \in \mathbb 
	Z_{\geq 0}$, the space $(F^{D_k}_{D_\omega}W_{D_k}\calO)(x, \ell)$ can be 
	described as the subspace of $(W_{D_\omega}\calO)(x,\ell)$ given by labeled 
	trees $f \colon y \to x$ with at most $\ell$ internal edges such that, 
	after cutting apart $f$ at the edges labeled $\infty$, each remaining 
	subtree of $f$ has $\leq k$ internal edges.  Taking the colimit as $k \to 
	\infty$ we get all labeled trees $f \colon x \to y$ with at most $\ell$ 
	internal edges.  That is, we get all of $(W_{D_\omega}\calO)(x, \ell)$.
\end{proof}

\begin{lemma}
  \label{lemma: FWkCO to Wk+1CO cofibration}
  Suppose $\calO \in \Op$ is cofibrant as an object of 
  $\fun^{\otimes}(\for_0,\topl)$ (see remark~\ref{remark: Sigma cofibrance and 
  for0}).  Then, for each $k\geq 0$ the natural map $\iota_{k} \colon 
  F^{D_{k-1}}_{D_k}W_{D_{k-1}}\calO \to W_{D_{k}}\calO$ is a cofibration in 
  $D_k\dash\Op$.
\end{lemma}

\begin{proof}
	Observe that if $f \colon x \to y$ is a tree in $D_k$ with at least one 
	internal edge, then $x \in D_{k-1}$, so $\iota_k(x)$ is an isomorphism.  By 
	lemma~\ref{lemma:when a cofibration of symmetric sequences is cofibration 
	of operads} we only need to show that $\iota_k$ is a cofibration of 
	$D_k$-symmetric sequences.

	In the case $k=0$,  if $x$ is a tree in $C = D_{-1}$ then $\iota_0(x,0)$ is 
	an isomorphism. If $x$ is in $D - C$ then $\iota(x,0)$ is the map 
	$\emptyset \to \calO(x)$ .  The assumption that $\calO$ is cofibrant as a 
	functor $\fun^{\otimes}(\for_0, \topl)$ implies in particular that each 
	$\calO(x)$ is $\aut(x)$-cofibrant.  That is, $\calO$ is cofibrant as a 
	$\for$-symmetric sequence, so $\iota_0$ is a cofibration of $D_0$-symmetric 
	sequences.  We now consider the case $k > 0$.  We follow \cite[lemma 
	5.4]{berger-moerdijk-w}.  For a young tree $z \in \for_k$ and $g \colon y 
	\to z$ a tree in $\for_{k}$, let $(W\times \calO)^+_k(g)$ be $W(g) \times 
	\calO(y)$ if $g$ has $\leq k-1$ internal edges.  Otherwise let $(W\times 
	\calO)^+_{k}(g) \subset W(g) \times \calO(y)$ be the set of $(t,\alpha)$ 
	such that $t(\epsilon) = 0$ or $t(\epsilon)=\infty$ for some $\epsilon \in 
	E(g)$ or $\alpha(j) = \id$ for some $j \in V(g)$.  	Using the techniques of 
	Berger and Moerdijk \cite[section 2]{berger-moerdijk-w} one can also show 
	that $\calO$ cofibrant as an object of $\fun^{\otimes}(\for_0, \topl)$ (see 
	remark~\ref{remark: Sigma cofibrance and for0}) then $(W\times 
	\calO)^+_k(g) \to W(g) \times \calO(\dom g)$ is an $\aut(g)$-cofibration, 
	where $\aut(g)$ is the automorphism group of $g$ as an object of the 
	category $\for_{k/z}$. Define the map $(W\times \calO)^+_{k}(g) \to 
	F^{D_{k-1}}_{D_k}W_{D_{k-1}}\calO(z)$ by collapsing any edge labeled $0$ 
	and deleting any vertex labeled with the identity.
	
	In the diagram below, $\pi_0 \for_{k/z}$ is the set of isomorphism classes 
	of forests $g \colon y \to z$ in $\for_k$.  The domain, $y$, of a forest $g 
	\colon y \to z$ is denoted $\dom g$.
	\begin{equation}
		\label{diagram:FWO-to-WO-pushout}
		\begin{tikzpicture}[
				numbered picture,
				node distance=8ex and 15em,
				display math nodes,
				text height=2ex,
				text depth=0.75ex,
				every node/.style={on grid,inner sep=1pt},
				mn/.style={inner sep=1pt, font=\scriptsize}
			]
			\node (ul) {\coprod_{[g] \in \pi_0 \for_{k/z}} 
			((W\times\calO)^+_{k}(g))_{\aut(g)} };
			\node[right=of ul] (ur) {(F^{D_{k-1}}_{D_k}W_{D_{k-1}}\calO)(z)};
			\node[below=of ul] (ll) {\coprod_{[g]\in\pi_0\for_{k/z}} (W(g) 
			\times \calO(\dom g))_{\aut(g)}};
			\node[right=of ll] (lr) {(W_{D_k}\calO)(z)};
			\path[->]
			(ul) edge (ur)
			(ul) edge (ll)
			(ur) edge (lr)
			(ll) edge (lr);
		\end{tikzpicture}
		\medskip
	\end{equation}
	The square in diagram~\eqref{diagram:FWO-to-WO-pushout} is a pushout.
	By remarks above, we know the map on the left in 
	\eqref{diagram:FWO-to-WO-pushout} is an $\aut(z)$-cofibration.  We conclude  
	$F^{D_{k-1}}_{D_k}W_{D_{k-1}}\calO \to W_{D_k}\calO$ is a cofibration of 
	$D_k$-symmetric sequences.
\end{proof}

\begin{lemma}
	\label{lemma: FWCO to WDO cofibration}
	If $C$ is a full subcategory of $D$ and $D$ is a full subcategory of 
	$\for$, and $\calO$ is an operad $\for \to \topl$ which is cofibrant as a 
	functor $\for_0 \to \topl$, then the natural map $F^C_DW_C\calO \to W_D 
	\calO$ is a cofibration of $D$-operads.
\end{lemma}
\begin{proof}
	By lemma~\ref{lemma:Fomega Womega is W} we have
	\[
		W_C\calO \cong F_C^{C_\omega} W_{C_\omega} \calO = F^C_{C_\omega} 
		W_{D_{-1}} \calO.
	\]
	The commutative diagram \eqref{diagram:sequence of forgetful functors} 
	shows that $F^C_DF^{C_\omega}_{C} = F^{D_\omega}_{D}F^{D_{-1}}_{D_\omega}$.  
	This gives the first equality below.
	\begin{align*}
		F_D^C F_C^{C_\omega} W_{C_\omega} \calO &=
		F^{D_\omega}_D
		F^{D_{-1}}_{D_\omega} W_{D_{-1}}\calO\\
		& \to F^{D_\omega}_D \colim_k F^{D_k}_{D_\omega} W_{D_k} \calO\\
		& = F^{D_\omega}_D W_{D_\omega}\calO\\
		& \cong W_D \calO.
	\end{align*}
	The equality in the third line above comes from lemma~\ref{lemma: colim 
	FWkO is WomegaO} and the fourth from lemma~\ref{lemma:Fomega Womega is W}.  
	The arrow above is a cofibration.  Indeed, by repeated application of 
	lemma~\ref{lemma: FWkCO to Wk+1CO cofibration} together with the fact that 
	each $F^{D_{k-1}}_{D_k}$ preserves cofibrations (remark~\ref{remark: F 
	preserves cofibrations}) we see that $F^{D_{k-1}}_{D_k} \cdots 
	F_{D_0}^{D_{-1}} W_{D_{-1}} \calO \to W_{D_k}\calO$ is a cofibration in 
	$D_k\dash\Op$.  Again, by remark~\ref{remark: F preserves cofibrations}
	$F^{D_k}_{D_\omega}$ and $F^{D_\omega}_D$ preserve cofibrations.  Thus the 
	arrow above is a cofibration in $D\dash\Op$.
\end{proof}

\subsection{Weak equivalence proof}
\label{section:weak-equivalence}
 This section contains a proof of
 \begin{theorem}\label{theorem:schinf-sch-weak-equivalence}
  The natural map $\schinf \to \sch$ is a weak equivalence of operads in 
  $\Op(\Coll{\topl})$.
\end{theorem}
 The idea of the proof is to consider the maps
$p_1 \colon \schinf(n,m) \to \schinf(n-1,m)$ and $p \colon\sch(n,m) \to
\sch(n-1,m)$ given by forgetting the $n^{th}$ disk.  By induction, we can
suppose $\schinf(n-1,m) \to \sch(n-1,m)$ is a weak equivalence.  We continue 
the induction by showing  that $p_1\inv(\alpha) \to p\inv(\alpha)$ is a weak 
equivalence for
every $\alpha \in \schinf(n-1,m)$.

To make the computation of $p_1\inv(\alpha)$ and $p\inv(\alpha)$ accessible, we 
will collapse the $n^{th}$ disk of $\alpha \in \sch(n,m)$ to a point.  Our goal 
in the next section is to make this precise.

\subsubsection{Defining \texorpdfstring{$\swcheesebullet(k,\ell|n,m)$}{ 
SCbullet(k,l|n,m)}}

\begin{figure}[ht]
\label{swcheesehat-def}
 \centering
 \includegraphics[width=\textwidth]{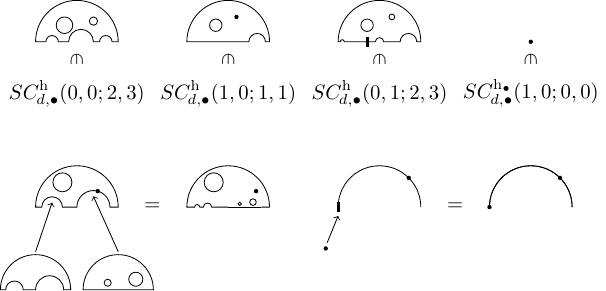}
\caption[The 3-colored operad $\swcheese^\h$]{The collapsed discs are denoted 
by dots and the collapsed half-discs by tick marks.  Collapsed discs are color 
$\e_\smallbullet$ input edges and collapsed half-discs are color 
$\h_\smallbullet$ input edges.  To keep the collapsed discs and half-discs from 
coinciding, we only allow one or the other in any composition.  Composition in 
$\swcheesebullet$ takes place only in the half-discs and collapsed half discs.  
The only composition we can do in a collapsed half-disk is given by plugging in 
a collapsed disk.  The result is a collapsed disk replacing the collapsed 
half-disk.}
\end{figure}

When we collapse the $n^{th}$ disk of $\alpha \in \swcheese^\h(n,m)$ to its 
center, we think of the result $\hat{\alpha}$ as living in a four-colored 
operad which we denote by $\swcheesebullet$.  We add the colors 
$\e_\smallbullet$ and $\h_\smallbullet$.  Let $K_\smallbullet = \{ 
\e_\smallbullet, \h_\smallbullet, \e, \h\}$ be the set of colors for this new 
operad.  The color $\e_\smallbullet$ stands for collapsed full disks.  It is 
convenient to also allow a collapsed half disk, which we color with 
$\h_\smallbullet$.  Let $(k,\ell|n,m)$ denote the $K_\smallbullet$-colored 
finite set with $k,\ell,n,$ and $m$ elements of color $\e_\smallbullet, 
\h_\smallbullet, \e,$ and $\h$ respectively.  Let
$D_\smallbullet$ denote the full sub category of $\rfor_{K_{\smallbullet}}$ 
with objects isomorphic to disjoint unions of the young trees
\begin{equation}\label{equation:smallbullet-connected-young-forests}
  (0,0|n,m) \to \{\h\} \hspace{0.6em} (1,0|n,m) \to \{\h\} \hspace{0.6em} 
  (0,1|n,m) \to \{\h\} \hspace{0.6em} (1,0|0,0) \to \{\h_\smallbullet\}.
\end{equation}
To define $\swcheesebullet\colon
D_\smallbullet \to \topl$ we need the notion of the geometric realization of 
$\beta\in \swcheese^\h(n, m)$.
\begin{definition}
\label{definition:geometric-realization}
 Given $\beta \in \swcheese^\h(n,m)$, let $\abs{\beta}$ be its \emph{geometric 
 realization}.  This is the subset of $\mathbb R^d$ given by deleting the open 
 discs and half-discs of $\beta$ from the closed unit half-disk.  More 
 precisely, if $\bar{D}^d_+$ is the closed unit half-disk in $\mathbb R^d$, 
 $\{(D^d_\e)_j\}_{j=1}^m$ are the open discs of $\beta$, and  
 $\{(D^d_\h)_i\}_{i=1}^{n}$ are the open half-discs of $\beta$ considered as 
 open discs in $\mathbb R^d$ whose center lies in $\mathbb R^{d-1}$, then
\[
 \abs{\beta} = \bar D^d_+ - \left(\left(\bigcup_{i=1}^{n} (D^d_\h)_i \right) 
 \cup \left(\bigcup_{\smash{j}=1}^m (D^d_\e)_j \right) \right).
\]
Let $\partial_\h \!\abs{\beta} \coloneqq \partial(\bar D^d_+ - (\cup_i 
(D^d_\h)_i))$ be the \emph{$\h$-colored boundary} of $\abs{\beta}$.  Let 
$\partial_{\rt}(\abs{\beta})$ be the upper hemisphere $S^{d-1}_+ \subset 
\partial \bar D^d_+$ and let $\partial_{i}\abs{\beta}$ be the upper hemisphere 
of $\partial(D^d_\h)_i$ for $1 \leq i \leq n$.
\end{definition}
Now we can set
\begin{align*}
\swcheesebullet^\h(0,0|n,m) &= \swcheese^\h(n,m) \\
\swcheesebullet^\h(1,0|n,m) &= \{(\alpha, q) \mid \alpha \in \swcheese^\h(n,m), 
q \in \abs{\alpha}\}  \\
\swcheesebullet^\h(0,1|n,m) &=\{ (\alpha, q) \mid \alpha \in \swcheese^\h(n,m), 
q \in \abs{\alpha}\cap \mathbb R^{d-1} \}\\
\swcheesebullet^{\h_\smallbullet}(k,\ell|n,m) &= \ast
\end{align*}
We think of the point $q \in \abs{\alpha}$ as a collapsed disk and the point $q 
\in \abs{\alpha}\cap \mathbb R^{d-1}$ as a collapsed half-disk.
Composition in $\swcheesebullet$ takes place in the half-discs and collapsed 
half-discs only. The un-collapsed discs play no part in composition.  However 
the collapsed half-discs and collapsed discs only play a part in composition 
when we plug a collapsed disk into a collapsed half-disk.  The result is a 
collapsed disk which happens to live on the boundary of the geometric 
realization.  See figure~\ref{swcheesehat-def}.

\begin{definition}\label{definition:e-leq-1-forest-categories}
Let
$C_\smallbullet$ denote the full sub category of
$D_\smallbullet$
with objects isomorphic to disjoint unions of the young forests from 
\eqref{equation:smallbullet-connected-young-forests} with $n \leq 1$. Recall 
from example~\ref{example:sch as w construction} that we used $D$ to denote  
the full subcategory of $\for_K$ whose young trees are isomorphic to the trees 
in \eqref{equation:smallbullet-connected-young-forests} with $k = \ell = 0$.  
Likewise $C$
is the full sub category of
$D$ given by trees isomorphic to disjoint unions of forests from 
\eqref{equation:smallbullet-connected-young-forests} with $k = \ell =0$ and $n 
\leq 1$.   We write $x \leq 1$ if $x \in C$ or $x \in C_\smallbullet$.
We have that $\schone$ is the restriction of $\swcheese$ to $C$, and $\sch$ is 
the restriction of $\swcheese$ to $D$.
Let $\swcheesebullet^{\h 1}$ denote the restriction of $\swcheesebullet$ to 
$C_\smallbullet$.
\end{definition}

Let $\schbullet$ and $\schonebullet$ denote the $\W$ construction applied to 
the four-colored operads $\swcheesebullet$ and $\swcheesebullet^{\h 1}$ via 
definition~\ref{definition:WvarO}.  Let $F$ denote Kan extension along 
$C_\smallbullet \to D_\smallbullet$, that is, $F = 
F^{C_\smallbullet}_{D_\smallbullet}$ using the notation from 
\ref{notation:FCD}.

\begin{definition}\label{definition:three-p's}
Define $p \colon
D_\smallbullet \to D$
by sending the $K_\smallbullet$-colored young forest $x \colon I_x \to J_x$ to 
the $K$-colored forest $px$ with
\[
	I_{px} = (I_x)_{\e, \h} \quad J_{px} =  (J_x)_{\e, \h},
\]
where for $K' \subset K_\smallbullet$ we set $I_{K'}= \col_I\inv(K')$, $\col_I 
\colon I \to K_\smallbullet$.   In 
\eqref{equation:smallbullet-connected-young-forests} we see that 
$x(I_{px})\subset J_{px}$ so that we can define $px$ as the restriction of $x$ 
to $I_{px}$.  Observe that $p(1,0|n,m) = (n,m)$.  If $f\colon y \to x$ is a 
forest, then $pf \colon py \to px$ is defined using $f$.  Since $f$ preserves 
the colorings $pf$ is indeed a forest from $py$ to $px$.  If $f$ is a morphism 
in
$C_\smallbullet$
then $pf$ is a morphism in
$C$.

If $\beta \in \swcheesebullet(z)$ for a $K_\smallbullet$-colored young forest 
$z$, then we get $p\beta \in \swcheese(pz)$.  To define $p\beta$ write $\beta = 
(\beta_j)_{j\in J_z}$ where $\beta_j \in \swcheesebullet(z\inv(j);j)$.  Each 
$\beta_j$ is of the form $(\gamma_j, q_j)$ with $q_j \in \abs{\gamma_j}$ or of 
the form $\gamma_j \in \swcheese(z\inv(j);j)$.  Set $p\beta = (\gamma_j)_{j \in 
J_{pz}}$.

If $t \in W(f)$ and $f \in
D_\smallbullet$, then $E(pf) \subset E(f)$ and $pt \in W(pf)$ is defined to be 
the pullback of $t \colon E(f) \to [0,\infty]$.

If $(f\colon y \to x, t \in W(f), \alpha \in \swcheesebullet(x))$ represents a 
point in $\schbullet(x)$ then $(pf, pt, p\alpha)$ represents a point in 
$\sch(x)$.  This defines the map $p \colon \schbullet(1,0|n,m) \to \sch(n,m)$.  
The restriction of $p$ to $F(\schonebullet)(1,0|n,m)$ factors through the 
inclusion $\schinf(n,m) \to \sch(n,m)$.  Let $p_1$ be the induced map 
$F(\schonebullet)(1,0|n,m) \to \schinf(n,m)$.
\end{definition}

Consider the commutative diagram of topological spaces where the horizontal 
arrows do not assemble to operad maps,
\[
\begin{tikzpicture}
  \matrix (m) [matrix of math nodes, row sep={2ex}, column sep={3ex}, text 
  height=1.5ex, text depth=0.25ex]
{ \schinf(n+1, m)  &   F(\schonebullet)(1,0| n, m) & \schinf(n, m)\\
 \sch(n+1, m) & \schbullet(1,0| n, m) & \sch(n, m).
 \\};
\path[->]
 (m-1-1) edge (m-2-1) edge node[above, map name] {$\sim$} (m-1-2)
 (m-1-2) edge (m-2-2) edge node[above, map name] {$p_1$} (m-1-3)
 (m-2-1) edge node[above, map name] {$\sim$} (m-2-2)
 (m-2-2) edge node[above, map name] {$p$} (m-2-3)
 (m-1-3) edge node[right, map name] {$\iota$} (m-2-3);
 \end{tikzpicture}
\]
The maps $p_1$ and $p$ delete the collapsed disk and, if necessary, a left over 
collapsed half-disk.
By induction on $n$ we assume the right vertical arrow is an equivalence. We 
will show that for each $\alpha\in \schinf(n, m)$ the inclusion 
$p_1^{-1}(\alpha)\rightarrow p^{-1}(\iota\alpha)$ is an equivalence. Then by 
the long exact sequence of homotopy groups we conclude that the middle vertical 
arrow is an equivalence. The top left and bottom right maps collapse the 
$n^{th}$ full disk.  One can show that these are equivalences. We conclude that 
the left vertical arrow is also an equivalence.  This will prove 
theorem~\ref{theorem:schinf-sch-weak-equivalence}.

\subsubsection{Computing \texorpdfstring{$p\inv(\iota\alpha)$ and 
$p_1\inv(\alpha)$}{p-1(iota alpha) and p1-1(alpha)}.}

We have shown that the proof rests on 
proposition~\ref{proposition:fiber-equivalence} below.  This section is 
dedicated to the proof of this proposition.
\begin{proposition}
\label{proposition:fiber-equivalence}
 Fix $\alpha \in \schinf(n,m)$.  The inclusion of the fiber $p_1\inv(\alpha)$ 
 into the fiber $p\inv(\iota \alpha)$ is a weak equivalence.
\end{proposition}



Combining the $W$ construction \eqref{equation:WO(z)-as-coend} and the left 
adjoint $C\dash\Op \to D\dash\Op$ (i.e. $\Op(\Collleqone{\topl}) \to 
\Op(\Coll{\topl})$)  we get
\[
	\schinf(n,m) = \Bigg(\coprod_{\substack{g \colon z\to y\\ f \colon y \to 
	(n,m)}}W(g) \times \swcheese(z) \Bigg)\modsim,
\]
where $y \leq 1$ (definition~\ref{definition:e-leq-1-forest-categories})
.  If $\alpha \in \schinf(n,m)$ is represented by $(f,g, t, \tilde\alpha)$ 
where $f \colon y \to (n,m)$, $g \colon z_\alpha \to y$, $t\in W(g)$, and 
$\tilde\alpha \in \swcheese(z)$, then $\iota \alpha \in \sch(n,m)$ is 
represented by $(fg, W_\infty(f)t,\tilde\alpha)$.  Let $T_\alpha = fg\colon 
z_\alpha \to (n,m)$ and $t_\alpha=W_\infty(f)t$.  Without loss of generality, 
we may assume $t_\alpha(i) > 0$ for every $i \in E(T_\alpha)$ and that 
$\tilde\alpha(j) \neq \id_{\swcheese}$ for any $j \in J_{z_\alpha}$.

\begin{definition}\label{definition:W-alpha}
Let $\trees(1,0|n,m)$ denote the over-category 
$(D_\smallbullet)_{/(1,0|n,m)}$. Similarly,
let $\trees(n,m) =
D_{/(n,m)}$.
Let $p \colon \trees(1,0|n,m) \to \trees(n,m)$ denote the functor induced by 
$p$ from definition~\ref{definition:three-p's}.

 Note that $T_\alpha \in \trees(n,m)$.  Let $(S, \nu) \in 
 \trees(1,0|n,m)_{/T_\alpha}$ where $S\colon x \to (1,0|n,m)$ is any 
 $K_\smallbullet$-colored tree and $\nu \colon z_\alpha \to px$ is a forest 
 such that $(pS)\nu = T_\alpha$. Define functors $W_\alpha \colon 
 \trees(1,0|n,m)_{/T_\alpha}^{op} \to \topl$ and $\swisscheese_\alpha \colon 
 \trees(1,0|n,m)_{/T_\alpha} \to \topl$ via the pullbacks
\begin{equation}\label{diagram:swisscheesealpha-and-Walpha}
	\begin{tikzpicture}[numbered picture]
	  \matrix (m) [matrix of math nodes, row sep=2ex, column sep=1em, text 
	  height=1.5ex, text depth=0.5ex]
	  {
		\swisscheese_\alpha(S) &             & \swcheesebullet(x) & & 
		  W_\alpha(S) & W(S)  & W(pS) \\
		\ast                   &|[xshift=-3em]|\swcheese(z) & \swcheese(px)      
		  & & \ast        &       & W(T_\alpha).
	  \\};
	  \foreach \source/\target in {1-1/1-3, 1-1/2-1, 2-1/2-2, 2-2/2-3, 1-3/2-3, 
	  1-5/1-6, 1-5/2-5, 2-5/2-7, 1-6/1-7, 1-7/2-7}
	  {\path[->] (m-\source) edge (m-\target);}
	  \foreach \source/\target/\pos/\maplabel in {1-3/2-3/right/p, 
	  2-1/2-2/above/\tilde\alpha, 2-2/2-3/above/\swcheese(\nu), 
  1-6/1-7/above/p, 1-7/2-7/right/W_\Sigma(\nu), 2-5/2-7/above/t_\alpha}
	  {\path (m-\source) edge node[map name, \pos] {$\maplabel$}
	  (m-\target);}
	\end{tikzpicture}
  \end{equation}
\end{definition}

We want to replace $\trees(1,0|n,m)$ by a much smaller category.  First we need 
the wedge operation on forests.
\begin{definition}
  \label{definition:wedging-forests}
  Let $f \colon x \to y$ be a $K_f$-colored forest and let $g \colon z \to w$ 
  be an $K_g$-colored forest for some finite sets $K_f, K_g$.  Let $\tau \colon 
  J_w \to J_x$ be any map.  Define $x \vee_\tau z$ to be the young $K_f \sqcup 
  K_g$-colored forest $(x,\tau, z) \colon I_x \sqcup J_w \sqcup I_z \to J_x 
  \sqcup J_z$ and define $y \vee_\tau w$ to be the young forest $(y,[f|x]
  \tau w) \colon I_y \sqcup I_w \to J_y$ (see 
  definition~\ref{definition:(f|x)}).  Finally, set $f\vee_\tau g \colon x 
  \vee_\tau z \to y \vee_\tau w$ to be the forest $(f,g,f,\tau)\colon I_y 
  \sqcup I_w \sqcup J_x \sqcup J_z \to I_x \sqcup J_w \sqcup I_z \sqcup J_y$.
\end{definition}

\begin{definition}
  \label{definition:trees-alpha}
  Let $\Gamma_0$ be the tree with with no internal vertices and a single input 
  vertex of color $\e_\smallbullet$.  Let $\Gamma_1$ be the tree with a single 
  internal vertex of color $\h_\smallbullet$ and a single input vertex of color 
  $\e_\smallbullet$.

  For any edge $i \in \exEdges(T_\alpha)_\h$ define $\nu(i) \colon T_\alpha \to 
  T_\alpha(i)$ to be the morphism in $\trees(n,m)$ which inserts a unary vertex 
  along $i$. Call this new vertex $i_v$.   Let $S_{i,k} = T_\alpha(i) 
  \vee_{i_v} \Gamma_k$.  For any internal vertex $j \in J_{z_\alpha}$ let 
  $S_{j,k} = T_\alpha \vee_{j} \Gamma_k$.  Note that $pS_{i,k} = T_\alpha(i)$ 
  and $pS_{j,k} = T_\alpha$.

  Let $\trees_\alpha$ be the full subcategory of $\trees(1,0|n,m)_{/T_\alpha}$ 
  given by the objects $S_{i,k}=(S_{i,k},\nu(i))$ and 
  $S_{j,k}=(S_{j,k},\id_{T_\alpha})$ where $i \in (I_{z_\alpha})_\h\sqcup\{ \rt 
  \}$, $j \in J_{z_\alpha}$ and $k \in \{ 0,1 \}$.
\end{definition}

\begin{figure}[htbp]
\label{figure:trees-over-Talpha}
\centering
\includegraphics[width=\textwidth]{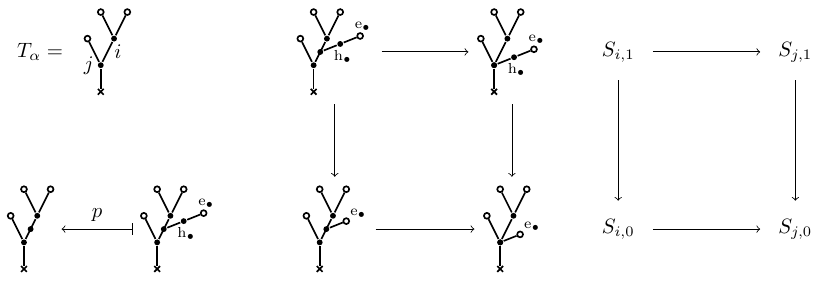}
 \caption[The category $\trees_{\alpha}$]{The edge $i$ and vertex $j$ of 
 $T_\alpha$ give a commutative square in $\trees_{\alpha}$.  The input vertices 
 are circles.  The output vertex ends in an $\mathsf x$.  The internal vertices 
 are filled dots.  The input and internal vertices of $\Gamma_{0}$ and 
 $\Gamma_{1}$ are labeled with their colors.  In addition, the image of 
 $S_{i,1}$ under the functor $p$ is shown.  This makes it clear that the map 
 $T_\alpha \to pS_{i,1}$ is given by inserting a single vertex.}
\end{figure}

\begin{remark}
The advantage of $T_\alpha$ is that it is easy to understand and computes the 
space $p\inv\alpha$ (lemma~\ref{lemma:p-inv-alpha-coend}).  There is a unique 
morphism $S_{\ell, 1}$ to $S_{\ell, 0}$ for every $\ell$ and unique morphisms 
$S_{i,k} \to S_{T_\alpha\inv(i),k}$ and $S_{i,k} \to S_{z_\alpha(i),k}$.   See 
figure \ref{figure:trees-over-Talpha} for an illustration.
\end{remark}

\begin{lemma}
\label{lemma:p-inv-alpha-coend}
 The fiber $p\inv(\alpha)$ is given by the coend
\[
 W_\alpha \otimes_{\trees_{\alpha}} \swisscheese_\alpha.
\]
\end{lemma}
\begin{proof}
Let $\gamma = [S, s, \tilde\gamma] \in \schbullet(1,0|n,m)$ where $S \colon x 
\to (1,0|n,m)$ is a forest in
$D_\smallbullet$, $s \in W(S)$, and $\tilde\gamma \in \swcheesebullet(x)$.  Let 
us assume that $\tilde\gamma(j)\neq \id$ for all $j \in J_x$ and $s(i) > 0$ for 
all $i \in E(S)$.  Observe that $p\gamma \in \sch(n,m)$ is given by 
$[pS,ps,p\tilde\gamma]$.  If $p\gamma=\alpha$ there must be some $\nu \colon 
T_\alpha \to pS$ in $\trees(n,m)$ such that $\swcheese(\nu)\tilde\alpha= 
p\tilde\gamma$ and $W(\nu)ps = t_\alpha$.  The condition $t_\alpha(i) >0$ for 
all $i\in E(T_\alpha)$ implies that $t_\alpha\neq W_\Sigma(\nu)(t')$ for any 
$t'$ and any $\nu$ which collapses any edges.  Moreover the condition 
$\tilde\gamma(j) \neq \id$ for all $j$ implies that $p\tilde\gamma(j) \neq \id$ 
for all $j\in J_x$ such that $x\inv(j)_{\e_\smallbullet, \h_\smallbullet}$ is 
not empty.  We conclude that either $\nu =\id$ or $\nu$ is the insertion of the 
unique unary (in $pS$, not in $S$) vertex $j$ such that 
$x\inv(j)_{\e_\smallbullet, \h_\smallbullet}$ is not empty.  In the former case 
we must have $S= S_{j,k}$ for some vertex $j \in J_x$ and some $k \in \{  0,1 
\}$  In the latter case we have $S=S_{i,k}$ for some edge $i$ of $T_\alpha$ and 
some $k$. This defines the map $p\inv(\alpha) \to W_\alpha 
\otimes_{\trees_\alpha} \swisscheese_\alpha$.  The map in the other direction 
is clear and the verification that they are inverses is left to the reader.
%
%
%
\end{proof}

In diagram~\eqref{diagram:swisscheese-alpha-double-cube} we have $\h$-colored 
edges $i_1, i_2$ of $T_\alpha$ with $z_\alpha(i_1) = j = T_\alpha\inv(i_2)$.  
Thus we get the commutative diagram on the left.  The image of this diagram 
under $\swisscheese_\alpha$ is shown on the right.
\begin{equation}\label{diagram:swisscheese-alpha-double-cube}
\begin{tikzpicture}[numbered picture]
	  \matrix (m) [matrix of math nodes, row sep={2ex}, column sep={0.9em},
		 text height=1.5ex, text depth=0.25ex]
	  {
	   S_{i_1,1} & S_{j,1} & S_{i_2, 1} & \abs{\id_{\h}} \cap \mathbb R^{d-1} & 
		  \abs{\tilde \alpha(j)}\cap  \mathbb R^{d-1} &  \abs{\id_{\h}}\cap 
		  \mathbb R^{d-1} \hspace{-2em}\\
	   S_{i_1,0} & S_{j,0} & S_{i_2,0} &  \abs{\id_{\h}} & \abs{\tilde 
	  \alpha(j)} &  \abs{\id_{\h}}
	  \\};
	  \foreach \source/\target in {1-1/1-2, 1-1/2-1, 2-1/2-2, 1-2/2-2, 1-3/1-2, 
	  1-3/2-3, 2-3/2-2, 1-4/1-5, 1-4/2-4, 2-4/2-5, 1-5/2-5, 1-6/1-5, 1-6/2-6, 
  2-6/2-5}
	  {\path[->] (m-\source) edge (m-\target);}
	  \foreach \source/\target/\pos/\maplabel in {}
	  {\path (m-\source) edge node[map name, \pos] {$\maplabel$}
	  (m-\target);}
   \end{tikzpicture}
\end{equation}
The geometric realization of the identity $\id_\h$ is just $S^{d-1}_+$, the top 
half of the $(d-1)$-sphere.  The input of $\tilde\alpha(j)$ corresponding to 
$i_1$ is a half disk and the map $\abs{\id_\h} \to \abs{\tilde\alpha(j)}$ 
corresponding to $S_{i_1, 0} \to S_{j,0}$ is just 
$\partial_{i_1}\abs{\tilde\alpha(j)} \to \abs{\tilde\alpha_j}$ (see 
definition~\ref{definition:geometric-realization}).  On the other hand the 
image of $S_{i_2, 0} \to S_{j,0}$ is the inclusion of the output boundary 
$\partial_\rt\abs{\tilde\alpha(j)} \to \abs{\tilde\alpha(j)}$.

\begin{definition}
  \label{definition:the-three-edges-of-S}
  Let $\epsilon_\smallbullet \in E(S_{\ell,1})$ be the unique internal edge of 
  color $\h_\smallbullet$.  If $i \in \exEdges(T_\alpha)$ 
  (\ref{definition:edges vertices of forest}), let $i_v$ denote the vertex 
  inserted by $\nu \colon T_\alpha \to pS_{i,k}$.  Let $i_\inverts$ and 
  $i_\outverts$ respectively denote the incoming and outgoing edges of $i_v$ 
  considered as internal edges of $S_{i,k}$.  For any object $S_{\ell,k}$ of 
  $\trees_\alpha$, let $E_\alpha(S_{\ell,k}) = \{ \epsilon_\smallbullet \}^k 
  \sqcup (\{ i_\inverts, i_\outverts \} \cap E(S_{\ell,k}))$.  This defines a 
  functor $E_\alpha \colon \trees_\alpha^{op} \to \mathrm{Set}$.
\end{definition}

The image under $W_\alpha$ of the square in 
diagram~\eqref{diagram:swisscheese-alpha-double-cube} is in 
diagram~\eqref{diagram:W-alpha-double-cube}.
\begin{equation}\label{diagram:W-alpha-double-cube}
\begin{tikzpicture}[numbered picture]
	  \matrix (m) [matrix of math nodes, row sep={2.2ex}, column sep={3em},
		 text height=1.7ex, text depth=0.25ex]
	  {
	   \lbrack0,\infty\rbrack^2 & \lbrack0,\infty\rbrack & 
		  \lbrack0,\infty\rbrack^2 \\
	   \lbrack0,\infty\rbrack   & \ast       & \lbrack0,\infty\rbrack
	  \\};
	  \foreach \source/\target in {1-2/1-1, 1-2/1-3, 2-1/1-1, 2-3/1-3, 2-2/1-2, 
	  2-2/2-1, 2-2/2-3}
	  {\path[->] (m-\source) edge (m-\target);}
	  \def\comma{,}
	  \foreach \source/\target/\pos/\maplabel in {1-2/1-1/above/(\id\comma 0), 
	  1-2/1-3/above/(\id\comma\infty), 2-1/1-1/left/(0\comma\id), 
  2-3/1-3/right/(0\comma\id), 2-2/1-2/right/0, 2-2/2-1/above/0, 
  2-2/2-3/above/\infty}
	  {\path (m-\source) edge node[map name, \pos] {$\maplabel$}
	  (m-\target);}
   \end{tikzpicture}
\end{equation}
More precisely,
\begin{equation}
  \label{equation:W-alpha-S-as-mapping-space}
  W_\alpha(S) = \{ s \colon E_\alpha(S)\to [0,\infty] \mid s(i_\inverts) + 
  s(i_\outverts) = t_\alpha(i) \},
\end{equation}
and $W_\alpha(S) \to W_\alpha(S')$ for a map $S' \to S$ in $\trees_\alpha$ is 
given by push forward of functions along the map of finite sets $E_\alpha(S) 
\hookrightarrow E_\alpha(S')$.  There is no condition on 
$s(\epsilon_\smallbullet)$, the length of the edge of color $\h_\smallbullet$.  
The isomorphism $W_\alpha(S_{i,1}) \to [0,\infty]^2$ sends $s$ to 
$(s(\epsilon_\smallbullet),r(s(i_\outverts), s(i_\inverts)))$ where
\[
	r(s_o, s_i) =  \frac{1-e^{-s_o}}{1 - e^{-s_i}},
\]
which lands in $[0,\infty]$ because $s_o + s_i = t_\alpha>0$. Note that $s_o=0$ 
if and only if $r(s_o,s_i) = 0$ and $s_o = t_\alpha$ if and only if $r(s_o, 
s_i)=\infty$.  Since the morphism $S_{i_1,1} \to S_{j,1}$ from 
diagram~\eqref{diagram:swisscheese-alpha-double-cube} collapses the edge 
$(i_1)_\outverts$ we get $W_\alpha(S_{j,1}) \cong \{ (r_\smallbullet,r) \in 
W_\alpha(S_{i_1,1}) \mid  r=0\}$.  In the same diagram, the morphism 
$S_{i_2,1}\to S_{j,1}$ collapses the edge $(i_2)_\inverts$, so we have 
$W_\alpha(S_{j,1}) \cong \{ (r_\smallbullet,r) \in W_\alpha(S_{i_2,1}) \mid 
r=\infty \}$.  The unique morphism $S_{i,1} \to S_{i,0}$ collapses the edge 
$i_{\smallbullet}$ so that $W_\alpha(S_{i,0}) \cong \{ (r_\smallbullet,r) \in 
W_\alpha(S_{i,1})\mid r_\smallbullet = 0 \}$.  The rest can be deduced from 
these cases.

\begin{lemma}
\label{lemma:W-alpha-hocolim}
 For any functor $F \colon \trees_\alpha \to \topl$ the coend $W_\alpha 
 \otimes_{\trees_{\alpha}} F$ is the homotopy colimit of $F$ over 
 $\trees_{\alpha}$.
\end{lemma}
\begin{proof}
It is clear from diagrams~\eqref{diagram:W-alpha-double-cube} and 
\eqref{diagram:swisscheese-alpha-double-cube} that $W_\alpha(S)$ is the geometric 
realization of the nerve of the under category of $S$ for each object $S \in 
\trees_\alpha$.  In addition the maps $W_\alpha(S) \to W_\alpha(S')$ for $S' 
\to S$ agree with the maps obtained from the nerves of under categories.
\end{proof}

\begin{lemma}
  \label{lemma:p-inv-alpha-is}
  We can explicitly compute $p\inv(\iota\alpha)$ as
  \[
	p\inv(\iota\alpha) \simeq \abs{\swcheese(T_\alpha)\tilde\alpha} \simeq 
	(S^{d-1})^{\vee n},
  \]
  where $\swcheese(T_\alpha)\tilde\alpha$ is the composition of all vertex 
  labels from $\iota\alpha$.
\end{lemma}
\begin{proof}
   Let $\trees_{\alpha,0}$ denote the full subcategory of $\trees_\alpha$ 
   consisting of objects $S_{j,0}$ and $S_{i,0}$ for internal vertices $j$ and 
   internal edges $i$.  This category is homotopy terminal, so by 
   lemma~\ref{lemma:W-alpha-hocolim} and lemma~\ref{lemma:p-inv-alpha-coend} we 
   have $p\inv(\iota\alpha) = \hocolim_{\trees_{\alpha,0}} 
   \swisscheese_\alpha$.  This is the same as the homotopy colimit of the 
   coequalizer diagram
	 \[
		\coprod_{i \in E(T_\alpha)} \abs{\id_\h} \rightrightarrows 
		\coprod_{j\in V(T_\alpha)} \abs{\tilde\alpha(j)},
	 \]
	 where one arrow is given by including into output parts of the boundaries 
	 of $\abs{\tilde\alpha(j)}$'s, and the other arrow is given by including 
	 into input boundaries.  These maps are cofibrations with disjoint images.  
	 Each space in the coequalizer diagram is cofibrant.  Thus the coequalizer 
	 diagram is already cofibrant as a functor $(\cdot \rightrightarrows \cdot) 
	 \to \topl$.  Thus we can compute the normal colimit.  It is clear that 
	 this is the same as composing the $\tilde\alpha(j)$'s via $T_\alpha$ then 
	 taking the realization of the result.  In addition $\abs{\beta}$ is 
	 equivalent to a wedge of $n$ spheres of dimension $d-1$ if $\beta \in 
	 \swcheese^{\h}(n,m)$.
\end{proof}

\begin{definition}
  \label{definition:W-alpha-1}
  Let $\trees_{\alpha,1}$ denote the full subcategory of $\trees_\alpha$ where 
  we discard the objects $S_{j,0}$ and $S_{i,0}$ for $j \in J_{z_\alpha}$ and 
  $i \in E(T_\alpha)$.  Define a functor $W_{\alpha,1} \colon 
  \trees_{\alpha,1}^{op}\to \topl$ by setting
  \[
	W_{\alpha,1}(S) = \{ s \colon E_\alpha(S) \to [0,\infty] \mid \sum_{i\in 
	E_\alpha(S)} s(i) = \infty \}.
  \]
\end{definition}

\begin{lemma}
  \label{lemma:p1-inv-alpha-as-a-coend}
  Suppose $t_\alpha < \infty$, and $n = 1$, then $p_1\inv(\alpha)$ is given by 
  the coend
  \[
	W_\alpha \otimes_{\trees_{\alpha,1}} \swisscheese_\alpha,
  \]
  where $\swisscheese_\alpha$ is the functor in 
  definition~\ref{definition:W-alpha} restricted to $\trees_{\alpha,1}$ and 
  $W_{\alpha,1}$ is defined in \ref{definition:W-alpha-1}.
\end{lemma}
\begin{proof}
  Let $\gamma\in F(\schonebullet)(1,m)$ such that $p_1(\gamma) = \alpha$.  Pick 
  a representative $(f,g,s, \tilde\gamma)$ where $f \colon y \to (1,0|1,m)$, 
  $y\leq 1$, $g\colon z \to y$, $s\in W(g)$ and 
  $\tilde\gamma\in\swcheesebullet(z)$.  Consider  $\iota\gamma \in 
  \schbullet(1,m)$, which is represented by $(fg,W_\infty(f)s, \tilde\gamma)$.  
  Recall that the condition $y \leq 1$  means that each connected component of 
  the young forest $y$ has at most one input whose color lives in $\{ \e, 
  \e_\smallbullet \}$.  This implies that $f$ has at least one internal edge $i 
  \in E(f)$.  Thus $W_\infty(f)s(i) =\infty$ when $i$ is viewed as an internal 
  edge in $fg$.

  We know $p\iota\gamma = \iota\alpha$, so $\iota\gamma$ is represented by some 
  triple $(S,s',\tilde\gamma)$ with $S \in \trees_\alpha$, $s' \in 
  W_\alpha(S)$, and $\tilde\gamma\in \swisscheese_\alpha(S)$.  The relations in 
  $\schbullet$ preserve edges of length $\infty$, so we must have $s'(i) = 
  \infty$ for some $i \in E(S)$.  We are assuming $t_\alpha (i) < \infty$ for 
  all $i \in E(T_\alpha)$, so the infinite edge in $S$ must be in 
  $E_\alpha(S)$.  This implies $s' \in W_{\alpha,1}(S)$.  Moreover we cannot 
  have such an infinite edge if $S = S_{j,0}$ for some vertex $j$ or 
  $S=S_{i,0}$ for some internal edge $i$.  Thus $S \in \trees_{\alpha,1}$.  
  This defines the map from $p\inv(\alpha)$ to the coend.  We leave the 
  remainder to the reader.
\end{proof}

\begin{lemma}
\label{lemma:W-alpha-1-hocolim}
 For any functor $F \colon \trees_{\alpha,1} \to \topl$ the coend $W_{\alpha,1} 
 \otimes_{\trees_{\alpha,1}} F$ is the homotopy colimit of $F$ over 
 $\trees_{\alpha,1}$.
\end{lemma}
\begin{proof}
  The argument here is similar to the proof of 
  lemma~\ref{lemma:W-alpha-hocolim}.
\end{proof}

\begin{corollary}
  \label{corollary:p1-inv-alpha-degree-1-is}
  If $t_\alpha < 0$ and $n =1$, then the fiber $p_1\inv(\alpha)$ is equivalent 
  to $\partial_\h\abs{\swcheese(T_\alpha)\tilde\alpha} \simeq S^{d-1}$.
\end{corollary}
\begin{proof}
  By the same argument as in lemma~\ref{lemma:p-inv-alpha-is}, 
  $\hocolim_{\trees_{\alpha,1}} \swisscheese_\alpha$ is equivalent to 
  $\colim_{\trees_{\alpha,1}} \swisscheese_\alpha$.  This is easily computed as 
  the $\h$-colored boundary of the composite of $\tilde \alpha$.
\end{proof}

\begin{proof}[Proof of {proposition~\ref{proposition:fiber-equivalence}}]
  Recall $\alpha$ is represented by $f \colon y \to (n,m)$, $y\leq 1$, $g 
  \colon z_\alpha \to y$, $t \in W(g)$ and $\tilde \alpha \in 
  \swcheese(z_\alpha)$.  By applying relations in $\schinf$ we may assume 
  $0<t<\infty$.  We may think of $(g, t)$ as representing an element of 
  $\schone(y)$ which we can write as $(\alpha(j))_{ j \in J_y}$.  If $\alpha(j) 
  \in \sch(n_j,m_j)$ then $n_j \leq 1$.  Clearly $p_1\inv(\alpha(j)) \simeq 
  p\inv(\alpha(j)) \simeq \ast$ when $n_j = 0$.  Since $t_{\alpha(j)} < \infty$ 
  we can use corollary~\ref{corollary:p1-inv-alpha-degree-1-is} to conclude 
  $p_1\inv(\alpha(j)) \simeq \partial_\h \abs{(\swcheese(g)(\alpha))(j)}$.  The 
  fiber $p_1\inv(\alpha)$ is equal to the colimit of the diagram
	\[
		\coprod_{i \in E(f)} \abs{1_\h} \rightrightarrows \coprod_{j \in V(f)} 
		p_1\inv(\alpha(j)),
	\]
	where one arrow is given by $\abs{1_\h} \simeq \partial_i\abs{\alpha(y(i))} 
	\to \partial_\h\abs{(\swcheese(g)(\alpha))(y(i))}$ and the other by 
	$\abs{1_\h}\simeq \partial_{\rt}\abs{\alpha(f(i))} \to 
	\partial_\h\abs{(\swcheese(g)(\alpha))(y(i))}$.  This colimit is clearly 
	$(S^{d-1})^{\vee n} \simeq p\inv(\iota\alpha)$.
\end{proof}

\begin{figure}[htbp]
 \centering
\label{figure:p-inv-alpha}
\includegraphics[width=\textwidth]{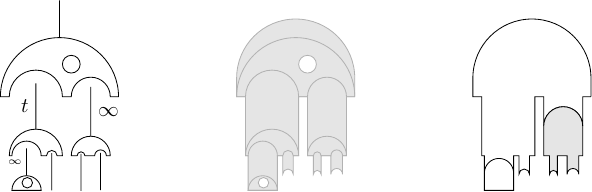}
\caption[Comparing $p_1\inv(\alpha)$ and $p\inv(\alpha)$]{On the left is
$\alpha \in \schinf(2,3)$.  In the middle is $p\inv(\alpha)$, and on the
right is $p_1\inv(\alpha)$.  Both $p\inv(\alpha)$ and $p\inv(\alpha)$ have
the homotopy type of a wedge of spheres, one for each disk in $\alpha$.}
\end{figure}

%
%
%
\bibliographystyle{plain}
\bibliography{jtbiblio}
%



\end{document}